\documentclass[10pt, a4paper]{article}
\usepackage[cp850]{inputenc}
\usepackage[all,poly,2cell,ps]{xy}

\usepackage{amssymb}
\usepackage{amsmath}

\usepackage{amsfonts}

\begin{document}

\author{G\'{e}rard Letac\thanks{Laboratoire de Statistique
et Probabilit\'es, Universit\'e Paul Sabatier, Toulouse, France.},   Mauro
Piccioni\thanks{Dipartimento di Matematica, Sapienza Universit\`{a} di Roma, 00185 Rome, Italia.}\ .} 

\title{Random continued fractions with beta-hypergeometric distribution}
\maketitle

\begin{abstract}In a recent paper (Asci \textit{et al.}, 2008) it has been shown
that certain random continued fractions have a density which is a
product of a beta density and a hypergeometric function $_{2}F_{1}$.
In the present paper we fully exploit a formula due to Thomae (1879)
in order to generalize substantially the class of random continuous
fractions with a density of the above form. This involves the design
of seven particular graphs. Infinite paths on them lead to random
continued fractions with an explicit  distribution. A careful study
about the set of five real parameters leading to a
beta-hypergeometric distribution is required, relying on almost
forgotten results mainly due to Felix Klein.


\end{abstract}

\section{\protect\bigskip Introduction} Recall that $_{2}F_{1}$ is the hypergeometric function defined as follows: the sequence of Pochhammer's symbols $\left\{
(t)_{n}\right\} _{n=0}^{\infty }$ is given for any real number $t$ by $%
(t)_{0}=1$ and $(t)_{n+1}=(t+n)(t)_{n}.$ For real numbers $p,q,r$ such
that $-r\notin \mathbb{N}=\{0,1,2,\ldots \}$ and for $0<x<1$  the number $_{2}F_{1}(p,q;r;x )$ is the  sum of the power
series 
\begin{equation}
\sum_{n=0}^{\infty }\frac{(p)_{n}(q)_{n}}{n!(r)_{n}}%
x ^{n}.  \label{hyperg}
\end{equation} Let $v=(a,b,p,q,r)^T$ be in $\mathbb{R}^5$, where $^T$ means transposition. Consider the function \begin{equation}h(a,b,p,q,r;x )=x ^{a-1}\left( 1-x \right) ^{b-1}\ 
_{2}F_{1}\left( p,q;r;x \right) ,\text{ }0<x <1,  \label{h}
\end{equation} and suppose that it is non negative and integrable on $(0,1).$ The distribution $BH(v)$ of $X\in (0,1)$ with density  proportional to $h(a,b,p,q,r;x )$
 is called a beta hypergeometric distribution. An important example is \begin{equation}\label{BHOne}BH(1,1,1,1,2)(dx)=\frac{6}{\pi^2 x}\log\frac{1}{1-x}\textbf{1}_{(0,1)}(x)dx.\end{equation}
If $c,d>0$ the probability measure  $$\beta^{(2)}_{c,d}(dw)=\frac{1}{B(c,d)}\frac{w^{c-1}}{(1+w)^{c+d}}\textbf{1}_{(0,\infty)}(w)dw$$ is called a beta distribution of the second kind. For instance if $U\sim \gamma_c$ and $V\sim \gamma_d$ are independent where the gamma distribution is  $$\gamma_c(du)=\frac{1}{\Gamma(c)}w^{c-1}e^{-u}\textbf{1}_{(0,\infty)}(w)dw$$ we have $U/V\sim \beta^{(2)}_{c,d}.$
Consider now two independent random variables $X\sim BH(v)$ and $W\sim \beta^{(2)}_{b',a'},$ where  
$a'=r-a$ and $b'=a+b-p.$ Also denote $p'=r+b-p-q,\ q'=b,\ r'=r+b-q.$ Therefore passing from $v=(a,b,p,q,r)^T$ to $v'=(a',b',p',q',r')^T$ is a linear automorphism  of $\mathbb{R}^5:$
$$v\mapsto v'=Mv$$
where $M$ is  a $5\times 5$  matrix (see (\ref{TM})). Similarly passing from $v=(a,b,p,q,r)^T$  to 
$$\Pi v=(a+b-p\, ,r-a)$$ is a linear operation. 
The present paper is based on the crucial fact that the random variable $(1+WX)^{-1}\in (0,1)$ has a density proportional to  $h(a',b',p',q',r';x ).$ In other terms we have the following 

\vspace{4mm}\noindent \textbf{Basic identity:} Let $X\sim BH(v)$ and $W\sim \beta^{(2)}_{\Pi v}$ be independent. Then 
\begin{equation}\label{BI}\frac{1}{1+WX}\sim BH(Mv).\end{equation} If $X\sim BH(v),$ $U\sim \gamma_{b'}$ and $V\sim \gamma_{a'}$ are independent an other presentation of the result is saying that $\frac{V}{V+UX}\sim BH(Mv).$  The  proof of the basic identity will be given by using the  formula obtained in 1879 by Thomae, which is an identity concerning the generalized hypergeometric function 
\begin{equation}\label{TFD}_3F_2(A,B,C;D,E;x)=\sum_{n=0}^{\infty }\frac{(A)_{n}(B)_{n}(C)_n}{n!(D)_{n}(E)_n}%
x ^{n}. \end{equation}Around this isolated basic identity, various considerations arise:

\vspace{4mm}\noindent \textsc{What are the beta hypergeometric distributions?} The exact knowledge of the acceptable set  of parameters $v=(a,b,p,q,r)^T$ such that the distribution $BH(v)$ exists needs a careful study of the positivity on $(0,1)$ of $\ _2F_1$. We rely for this problem on a 1890 paper by Felix Klein studying the number of zeros of $\ _2F_1$. Section 2 is devoted to a detailed description of a set $P$ of parameters $v=(a,b,p,q,r)^T$ such that the distribution $BH(v)$ exists. Actually the  $BH(v)$'s exist on a  set of $v$'s  which is slighly larger than $P;$  this larger set has an involved description and is not really useful to our purposes. 

\vspace{4mm}\noindent\textsc{Identifiability of the beta hypergeometric distributions.}
Things are complicated by the fact that the same beta hypergeometric distribution can be generically represented in four ways, first because of the symmetry $(p,q)$ and more importantly by the existence of the Euler formula for $\ _2F_1:$   
\begin{equation}
_{2}F_{1}\left( p,q;r;x \right) =\left( 1-x \right) ^{r-p-q}\ 
_{2}F_{1}\left( r-p,r-q;r;x \right) .  \label{euler}
\end{equation}
(see Rainville (1960) page 60).  In other terms consider the two $5\times 5$ matrices $T$ and $S$ defined respectively by 
  $(a,b,q,p,r)^T=T(a,b,p,q,r)^T$  and $(a,b+r-p-q,r-p,r-q,r)=S(a,b,p,q,r)^T$. The symmetry between $p$ and $q$ implies $BH(v)=BH(Tv)$ and the Euler formula implies $BH(v)=BH(Sv).$
  Actually the group generated by $T$ and $S$ has four elements $\{I_5,T,S,TS\}$ and is isomorphic to $(\mathbb{Z}/ 2\mathbb{Z})^2.$ Therefore the same beta hypergeometric distribution $BH(v)$ of $X$ can be generically represented in four ways by $v,$ $Tv,$ $Sv$ and $TSv$. On the other hand when $v=(a,b,p,q,r)^T$ we have
$$(a+b-p,r-a)=\Pi v=\Pi TS v,\ \ (a+b-q,r-a)=\Pi Tv=\Pi S v.$$
This implies that given $v$ there are only two ways, and not four, to use the basic indentity:  

\vspace{4mm}\noindent \textbf{Basic identity, revisited:} Let $X\sim BH(v)$ and $W\sim \beta^{(2)}_{\Pi v}=\beta^{(2)}_{\Pi TS v}$ be independent. Then 
$$\frac{1}{1+WX}\sim BH(Mv)=BH(MTSv).$$
If $X\sim BH(v)$ and $W\sim \beta^{(2)}_{\Pi T v}=\beta^{(2)}_{\Pi S v}$ are independent, then 
$$\frac{1}{1+WX}\sim BH(MTv)=BH(MSv).$$
The reader can immediately  verify that $MTS=TM$ and $MS=TMT.$

\vspace{4mm}\noindent\textsc{Explicit distributions of some random continued fractions.} 
The iteration of this basic identity leads to various types of random continued fractions which are described in Section 4. They  are splitted in seven categories corresponding to the seven partitions of the number 5, which are 
\begin{equation}\label{PART}5, 4+1,3+2, 3+1+1, 2+2+1, 2+1+1+1,1+1+1+1+1.\end{equation}
For instance as an application of Theorem 4.2 we see that since $v=(1,1,1,1,2)^T$ is an eigenvector of $M$ for the eigenvalue 1 and that $\Pi v=(1,1).$ Therefore if the $W_n$'s with $n=1,2\ldots$ are iid such that  $\Pr(0<W_n<w)=\frac{w}{1+w}$ (that is if $W_n\sim \beta^{(2)}_{1,1}$) we get that $$X=\frac{1}{1+\frac{W_1}{1+\frac{W_2}{...}}}$$
has distribution $BH(v)$ given by (\ref{BHOne}). This example belongs to the partition 5.

To sum up, Section 2 is a thorough study of $BH$ distributions, Section 3 proves the basic identity and introduces a crucial reparameterization of the BH laws for which in particular symmetries due to the Euler formula appear clearly. Section 4 applies the previous material to random continued fractions after introducing the graphs associated to a $BH$ law. This section contains  Theorem 4.2, the main result of the paper. 
Finally Section 5 is a detailed study of the cycles of the graphs defined on Section 4, which makes simple the application of Theorems 4.2 and 4.3. 

\section{Beta-hypergeometric distributions}
   This  long section (which can be skipped in first reading) describes the class of  probability densities on $(0,1)$ with five real parameters $(a,b,p,q,r)$ called beta hypergeometric densities. They have the form $x\mapsto C h(a,b,p,q,r;x )$ where $h$ is defined by (\ref{h}) 
   and  where $C$ is a suitable constant. 
To achieve this description of the beta hypergeometric densities we have to answer to the following questions
\begin{enumerate}
\item \textsc{Positivity problem}: when $_{2}F_{1}\left( p,q;r;x \right)\geq 0$ on $(0,1)?$ The set of acceptable $(p,q,r)$ will be denoted by $\mathbb{P}.$
\item \textsc{Integrability problem}: if $h=h(a,b,p,q,r;x )\geq 0$ when do we have $\int_0^1hdx<\infty?$ We shall define at the end of  Subsection 2.2 below a set $P$  of acceptable $(a,b, p,q,r)$ such that this is fulfilled, with $(p,q,r)$ is in $\mathbb{P}$. 
\item \textsc{Identifiability problem}: if $h(a^{\ast },b^{\ast },p^{\ast },q^{\ast },r^{\ast };x
)$ and $h(a,b,p,q,r;x )$ are proportional what are the relations between $(a^{\ast },b^{\ast },p^{\ast },q^{\ast },r^{\ast})$ and $(a,b,p,q,r)?$ 
\end{enumerate}
For doing so we recall a few classical facts about the hypergeometric function.
It is clear that if $p$ (or $q$) is a negative integer, then the power
series (\ref{hyperg}) happens to be a polynomial. In particular, if $p=0$,
then 

\begin{equation}
_{2}F_{1}(0,q;r;x )=1,  \label{bet}
\end{equation}%
for any value of $q$ and $r$, whereas if $p=-1$ then 

\begin{equation}
_{2}F_{1}(-1,q;r;x )=1-\frac{q}{r}x ,  \label{quasi}
\end{equation}%
meaning that for $q=sr$, the function (\ref{quasi}) is always equal to $%
1-sx $. Next recall that, provided $-t\notin \mathbb{N}$,  we can express 
$(t)_{n}$ as $\Gamma \left( t+n\right) /\Gamma \left( t\right) .$  Such a formula gives meaning to $\Gamma(t)$ when $t<0$ and  $t\notin-\mathbb{N}.$ Since actually $z\mapsto 1/\Gamma(z)$ is an entire function we write $1/\Gamma(-n)=0$ when $n\in \mathbb{N}.$ 

When $
p $ and $q$ are not negative integers we can apply Stirling's approximation
to the Gamma function to evaluate the order of the general term of the power
series (\ref{hyperg}). It turns out that it is equivalent (up to a
multiplicative constant) to $n^{-1-\left( r-p-q\right) }$. This implies that
its radius of convergence is equal to $1$, therefore the series is certainly
well defined for $x \in $ $\left[ 0,1\right) .$ It is convergent in $%
x =1$ if and only if $r-p-q>0$,  in which case (see Rainville (1960) page 49):

\begin{equation}\label{sophie}  \lim_{x\uparrow 1}\ _{2}F_{1}\left( p,q;r;x \right)=\frac{\Gamma(r)\Gamma(r-p-q)}{\Gamma(r-p)\Gamma(r-q)}\end{equation}
Note that this implies that this limit is zero if either $p-r\in \mathbb{N}$ or $q-r\in \mathbb{N}.$

F

\subsection{Positivity} The well known fact that $_{2}F_{1}(p,q;r;x)=z(x)$ is an analytic solution of the second order differential equation 
\begin{equation}
x\left( 1-x\right) z^{\prime \prime }+[r-\left( p+q+1\right) x]z^{\prime
}-pqz=0  \label{euler2}
\end{equation}%
in the unit disk  shows that the zeros of $_{2}F_{1}(p,q;r;x)$ on (0,1) are simple (for if $0=z(x_0)=z'(x_0)$ we easily see by induction using (\ref{euler2}) that $z^{(k)}(x_0)=0$ for all $k\in \mathbb{N}$ and thus $z\equiv 0$ which contradicts $z(0)=1).$ 
For this reason $_{2}F_{1}(p,q;r;x)\geq 0$ for $x\in (0,1)$ if and only if $_{2}F_{1}$ has no zeros in $(0,1).$ 
Denote \begin{equation}\label{P}\mathbb{P}=\{(p,q,r)\  ;\ \in \mathbb{R}^3\ ; r\notin -\mathbb{N}, \ F(p,q;r;x)>0 \ \forall x\in (0,1)\}\end{equation}
The aim of this subsection is to describe $\mathbb{P}.$

The determination of the number of zeros of $_{2}F_{1}$  is a difficult question which has been investigated in a number of papers, including Klein (1890), Hurwitz (1891), Van Vleck (1902).  To this aim  Klein introduces 
the function $x\mapsto E(x)$ on the real line defined by $E(x)=0$ if $x\leq 1$ and $E(x)=N$ if $N$ is the positive integer such that $N< x\leq N+1.$  Klein finally  introduces the non negative integer 
$$X=X(p,q,r)=E(\frac{1}{2}(|p-q|-|r-1|-|r-p-q|+1))$$ and  proves the following theorem

\vspace{4mm}\noindent\textbf{Theorem 2.1.} Suppose that $r\notin -\mathbb{N}.$ Then the number of zeros $ Z$ of $x\mapsto F(p,q;r;x)$ in $(0,1)$ is either $X$ or $X+1.$ 

\vspace{4mm}\noindent
This theorem leads to the description of $\mathbb{P}$ for which  we introduce the following notation: for $x\in \mathbb{R}$ the number $s(x)\in \{-1,0,1\}$
is $\mathrm{sign}\; 1/\Gamma(x)$ with the convention $\mathrm{sign}(0)=0.$ Therefore $s(x)=1$ for $x>0$ or $-2k-2<x<-2k-1$ with $k\in \mathbb{N}$,  $s(x)=-1$ for  $-2k-1<x<-2k$ with $k\in \mathbb{N}$, $s(x)=0$ for $x\in -\mathbb{N}.$ We need also to consider the set

\begin{equation}\label{S}S=\{(x,y,z)\in \mathbb{R}^3\ ; \ s(x)s(y)s(y)\geq 0\}\end{equation}
Therefore the part of $S$ which is in the octant $x,y,z\leq 0$  is a union of unit cubes; in the octant $x,y\leq 0$, $z\geq 0$ the set $S$ is a union of columns with unit square section; in the octant $x\leq 0$, $y,z\geq 0$ the set $S$ is the union of slices of unit height; in the octant $x,y,z\geq 0$ $S$ is the octant itself.

\vspace{4mm}\noindent\textbf{{Theorem 2.2}.} Suppose that $r\notin -\mathbb{N}.$ Then   $(p,q,r)\in \mathbb{P}$ if and only if 

\begin{enumerate}
\item either $r\geq 1$,  $r-p-q\leq 0$ and $p,q\geq 0;$ 
\item or $r\geq 1,$ $r-p-q\geq 0$ and $r-p,r-q\geq 0;$
\item or $r<1,$ $r-p-q\leq 0,$  $r-p,r-q\leq 1$ and  $(p,q,r)$ is in $S;$ 
\item or $r<1,$ $r-p-q\geq 0,$ $p,q\leq 1$ and  $(r-p,r-q,r)$ is in $S.$ 
\end{enumerate}

\vspace{4mm}\noindent\textbf{Remark.} The most obvious case of positivity is $p,q,r>0.$ It can be checked that it fulfills the sufficient conditions above to be in $\mathbb{P}.$ The same holds for $r-p,r-q,r>0.$ 

\vspace{4mm}\noindent\textbf{Proof.} From Klein's Theorem 2.1 a necessary condition for $(p,q,r)\in \mathbb{P}$ is $X(p,q,r)=0.$ We now consider three cases.

 \vspace{4mm}\noindent \textsc{The case $r\geq 1$.}$\Leftarrow:$ 
Trivially if $p,q\geq 0$ and $r\geq 1$ we have $F(p,q;r;x)\geq 0$ for all $x\in (0,1).$ If $r-p,r-q\geq 0$ Euler formula (\ref{euler}) provides the same result. $\Rightarrow:$ From Theorem 2.1 we have $X(p,q,r)=0$ , namely $\frac{1}{2}(|p-q|-|r-1|-|r-p-q|+1))\leq 1$ or \begin{equation}\label{antoine}|p-q|\leq r+|r-p-q|\end{equation} Now we can assume $p\geq q.$ If $r-p-q\geq 0$ (\ref{antoine}) implies $r-p,r-q\geq 0.$ If $r-p-q\leq 0$ (\ref{antoine}) implies $p,q\geq 0.$

\vspace{4mm}\noindent \textsc{The case $r< 1$ and $r-p-q\geq 0.$}$\Leftarrow:$ If $p,q\leq 1$ this implies $X=0$ since \begin{equation}\label{mathilde}X(p,q,r)=0\Leftrightarrow
\frac{1}{2}(|p-q|-|r-1|-|r-p-q|+1))=\max (p,q)\leq 1.\end{equation}  If $r-p-q> 0$ then from (\ref{sophie}) the limit of $F(p,q;r;x)$ as $x$ tends to 1 is $\frac{\Gamma(r)\Gamma(r-p-q)}{\Gamma(r-p)\Gamma(r-q)}$ If $r-p$ and $r-q$ are not in $-\mathbb{N}$ this number $\frac{\Gamma(r)\Gamma(r-p-q)}{\Gamma(r-p)\Gamma(r-q)}$ is $>0$ since $(r-p,r-q,r)$ is in $S.$  Therefore $Z$ cannot be $X+1=1.$ 

If $r-p$ or $r-q$ is in $-\mathbb{N}$, things are more complicated since in this case $F(p,q;r;1)=0$ by (\ref{sophie}). We introduce $Q(x)=x^{1-r}\ _2F_1(p+1-r,q+1-r;2-r;x)$ which is an other solution of (\ref{euler2}) which is linearly independent of $ \ _2F_1(p,q;r;x)$ Applying the part $r\geq 1$ of this theorem to the function  $_2F_1(p+1-r,q+1-r;2-r;x)$ shows that $Q(x)>0$ on $(0,1)$ (for this, observe that $X(p,q,r)=X(p+1-r,q+1-r,2-r)$ and thus $X(p+1-r,q+1-r,2-r)=0).$ Now assume that $Z=1,$ denote by $x_0$ the zero of $ \ _2F_1(p,q;r;x)$ in $(0,1)$ and apply the Sturm Liouville theorem to the interval $(x_0,1)$ and to the pair of solutions $ \ _2F_1(p,q;r;x),Q(x)$ of (\ref{euler2}). Thus $Q$ should have a zero in $(x_0,1)$, a contradiction. Therefore $Z=X=0.$ 
The remaining case $r-p-q= 0$  is obtained by passing to the limit.

$\Rightarrow:$ Since $(p,q,r)\in \mathbb{P}$ we have $X(p,q,r)=0$ and $p$ and $q\leq 1$ from (\ref{mathilde}). If $r-p-q> 0$ we use the fact that the limit in $1$ of $ \ _2F_1(p,q;r;x)$ as $x$ tends to 1 is $\frac{\Gamma(r)\Gamma(r-p-q)}{\Gamma(r-p)\Gamma(r-q)}$  and is $\geq 0.$ This implies $(r-p,r-q,r)$ is in $S.$ Again, the remaining case $r-p-q= 0$  is obtained by passing to the limit. 

\vspace{4mm}\noindent \textsc{The case $r< 1$ and $r-p-q\leq  0.$} This case is reduced to the preceding one by the Euler formula (\ref{euler}). $\square$

\subsection{Integrability}

\vspace{4mm}\noindent \textbf{Proposition 2.3.} Let $(a,b,p,q,r)\in\mathbb{R}^5$ such that $r\notin -\mathbb{N}.$ Assume that $h(a,b,p,q,r;x)\geq 0$ on $(0,1),$ that is $(p,q,r) \in \mathbb{P}$ defined by (\ref{P}).  The condition \begin{equation}\label{ABP}a,b,r+b-p-q>0\end{equation}  is sufficient for having 
$I=\int_0^1 h(a,b,p,q,r;x)dx$  finite.  Under these circumstances 
\begin{equation}
I =\frac{\Gamma \left( a\right) \Gamma
\left( b\right) }{\Gamma \left( a+b\right) }\sum_{n=0}^{\infty}\frac{\left( p\right)
_{n}\left( q\right) _{n}\left( a\right) _{n}}{n!\left( r\right) _{n}\left(
a+b\right) _{n}}  \label{term} 
=B(a,b)\ _{3}F_{2}\left( p,q,a;r,a+b;1\right) .  
\end{equation}
where $B\left( a,b\right) =\frac{\Gamma \left( a\right) \Gamma \left(
b\right) }{\Gamma \left( a+b\right) }$, whereas the second factor $_{3}F_{2}$
is a generalized hypergeometric function, evaluated at $x =1$.  Conversely if 
 
 \begin{enumerate}\item either $r-p-q\geq 0,$  $r-p\notin -\mathbb{N}$  and  $r-q\notin -\mathbb{N}$
 \item or $r-p-q\leq 0$, $p\notin -\mathbb{N}$  and  $q\notin -\mathbb{N}$\end{enumerate}
 (\ref{ABP}) is also necessary to have  $I$ finite.

\vspace{4mm}\noindent \textbf{Proof.}  If $a>0$ trivially $\int_0^{1/2}hdx$ is finite. If $r-p-q>0$ then $ \ _2F_1(p,q;r,x)$ has a finite limit in 1 and $b>0$ implies that $\int_{1/2}^1hdx$ is finite. If $r-p-q<0$ by (\ref{euler2}) and the fact that $r+b-p-q>0$ we similarly get that $\int_{1/2}^1hdx$ is finite. If $r-p-q=0$ note that the general term of the defining series for $h(a,b,p,q,r;x )$ has a constant sign. Therefore monotone convergence is applicable and the sum is $B(a,b) \ _{3}F_{2}\left( p,q,a;p+q,a+b;1\right).$ But the series $ _{3}F_{2}\left( p,q,a;p+q,a+b;1\right)$ converges because $p+q+a+b-p-q-a=b>0.$ This reasoning shows also that (\ref{term}) holds. Therefore the direct part is proved.

Conversely  if $r-p-q>0$ if none of  $r-p$ and $r-q$ is in $-\mathbb{N}$ we can claim that the limit of $ \ _2F_1(p,q;r;x)$ for $x\rightarrow 1$ exists and is not zero. Thus if furthermore $\int_0^1 h(a,b,p,q,r;x)dx$ is finite clearly $a$ and $b$ are positive, and $r+b-p-q>0$ a fortiori. The case $r-p-q<0$ is similar by the Euler formula (\ref{euler}). For $r-p-q=0$ we use the fact that  if neither $p$  nor $q$ are in $-\mathbb{N}$ we have 
\begin{equation}\label{nora} \ _2F_1(p,q;p+q;x)\sim_{x\rightarrow 1} \frac{1}{B(p,q)}\log\frac{1}{1-x}.\end{equation}
 To see this apply a classical result about the asymptotic behavior of a power series near the circle of convergence (see Titchmarsh (1939, 7.5))  to $$a_n=\frac{\Gamma(p+n)\Gamma(q+n)}{\Gamma(p+q+n)}\sim b_n=\frac{1}{n}. $$
 
 Note that if $(p,q,p+q)$ is in $\mathbb{P}$ defined by (\ref{P}) Theorem 2.2 implies that for $p+q\geq 1$ then $p$ and $q$ are positive and for $p+q<1$ we must have $(p,q,p+q)$  in $S$; thus $B(p,q)>0$ even when $p$ or $q$ are negative. Clearly (\ref{nora}) implies that $h$ is integrable  only if $a,b>0.$ $\square$

For the sake of completness in the following proposition we consider the cases where conditions 1 and 2 of the previous proposition are violated.

\vspace{4mm}\noindent \textbf{Proposition 2.4.} Let $(a,b,p,q,r)\in\mathbb{R}^5$ such that $r\notin -\mathbb{N}.$ Assume that $h=h(a,b,p,q,r;x)\geq 0$ on $(0,1)$ or equivalently that $(p,q,r)$ is in $\mathbb{P}.$ \begin{enumerate}
\item If $r-p-q\geq 0$ and $r-p$ or $r-q$ is $-n$ with $n\in \mathbb{N}$ then $\int_0^1hdx<\infty$ if and only if $a,r+b-p-q>0$ (but $b\leq 0$ is allowed)   

\item If $r-p-q\leq 0$ and $p$ or $q$ is $-n$ with $n\in \mathbb{N}$ then $\int_0^1hdx<\infty$ if and only if $a,b>0$ (but $r+b-p-q\leq 0$ is allowed)
 
 \end{enumerate}

\vspace{4mm}\noindent \textbf{Proof.} Parts 1 and 2 are obvious. To show Part 4, observe that if $r-p-q\leq 0$ and  $h\geq 0$ on $(0,1),$  Theorem 2.2 implies 
$r<1$ since $p=-n$ (say) is negative. The inequalities of Theorem 2.2 give $r-p=r+n\leq 1$ and $r-p-q=r+n-q\leq 0.$ For simplicity we write $s=r+n\leq 1$. This implies $s\leq q.$ Since $r=1-n$ is forbidden  by $r\notin -\mathbb{N},$ we have $s<1.$ Now 
$$ \ _2F_1(-n,q;r;1)=\frac{n!}{(1-s)_n}\sum_{k=0}^n\frac{(q)_k}{k!}\frac{(1-s)_{n-k}}{(n-k)!}=\frac{(q-s+1)_n}{(1-s)_n}>0.$$
Thus $\int_0^1hdx<\infty$ if and only if $a,b>0.$ Part 3 is proved from part 4 by using the Euler equality (\ref{euler}). 

\vspace{4mm}\noindent \textbf{Remark.} 
 From the two previous results the cases where $\int_0^1hdx<\infty$ does not imply $\min (b,r+b-p-q)>0$ are truly exceptional. An example is obtained from Proposition 2.4 part 1, with $r=p$ and $q=-2$; in this case 
 $\int_0^1hdx<\infty$ if and only if $a,b+2>0.$ Thus for instance $b=-1$ is allowed. For this reason we coin the following definition, using (\ref{TFD}):

 \vspace{4mm}\noindent \textbf{Definition.} The beta hyperbolic distribution $$BH(a,b,p,q,r)(dx)=
\frac{x^{a-1}(1-x)^{b-1}\ _2F_1(p,q;r;x) }{B(a,b)\ _3F_2(a,p,q;r;a+b;1)}\textbf{1}_{(0,1)}(x)dx$$ is defined if and  only if
\begin{enumerate}\item $r\notin -\mathbb{N};$
\item $(p,q,r)$ is in $\mathbb{P}$ described in  Theorem 2.2; 
\item $a>0$, $b>0$ and $r+b-p-q>0.$ \end{enumerate} In the sequel we denote by $P$ the set of $(a,b,p,q,r)$ satisfying the three conditions above.

\vspace{4mm}\noindent \textbf{Remark.} It is immediately verified
that the set $P$ is invariant under both $T$ and $S$. Recalling the
remark following Theorem 2.2 it is also immediately obtained that
the posititity of all the components of the vector $v=(a,b,p,q,r)^T$
(or $Sv$) together with the condition $r+b-p-q>0$ (or $b>0$) ensure
that $v \in P$.

\subsection{Identifiability}
 This subsection adresses to the problem of the identifiability, since we already know that \begin{eqnarray*}&&BH(a,b,p,q,r)=BH(a,b,q,p,r)=\\&&BH(a,r+b-p-q,r-p,r-q,r) =BH(a,r+b-p-q,r-q,r-p,r)\end{eqnarray*}from the obvious symmetry in $(p,q)$ and from Euler identity. Therefore BH distributions could have four different representations. On the other hand, formulas (\ref{bet}) and (\ref{quasi}) show that the number of representations of the same BH distribution can even be infinite.  We are interested in deciding when a function $h(a,b,p,q,r;x )$ defined by (\ref{h})
can be represented with different values of the parameters $a,b,p,q$ and $r.$
 The following result says essentially that,
aside from the symmetry of $_{2}F_{1}$ in $p$ and $q$ and the Euler's
identity (\ref{euler}), the only lack of a unique representation is due to the relations (\ref{bet}) and (\ref{quasi}). The theorem does not use the results of Sections 2.1 and 2.2.

\vspace{4mm}\noindent\textbf{Theorem 2.5.} Suppose that 
\begin{equation*}
h(a^{\ast },b^{\ast },p^{\ast },q^{\ast },r^{\ast };x
)=Ch(a,b,p,q,r;x ),0<x <1.
\end{equation*}%
Then $C=1$, $a^{\ast }=a$ and

\begin{description}
\item[a] Either $r^{\ast }=r$, and 
\begin{enumerate}
\item either $b^{\ast }=b,\left\{ p^{\ast },q^{\ast }\right\} =\left\{
p,q\right\} $,
\item or $b^{\ast }=b+r-p-q,\left\{ p^{\ast },q^{\ast }\right\} =\left\{
r-p,r-q\right\} $;
\end{enumerate}
\item[b] Or the function has the form%
\begin{equation}
h(a,b,p,q,r;x )=x ^{a-1}\left( 1-x \right) ^{u-1}(1-vx ).
\label{except}
\end{equation}

\begin{enumerate}
\item If $v=0$ either one between $p$ and $q$ is zero and $u=b$, or one
between $p$ and $q$ is equal to $r$, say $p$, and $u=b-q$ (the same obviously holds for 
$b^{\ast }, p^{\ast }, q^{\ast }$ and $r^{\ast });$

\item If $v\neq 0,1$ (which can be assumed without lack of generality)
either one between $p$ and $q$ is equal to $-1$, say $p$, $u=b$ and $v=\frac{%
q}{r}$, or one between $p$ and $q$ is equal to $r+1$, say $p$, $u=b-q-1$ and 
$v=\frac{r-q}{r}$(the same obviously holds for 
$b^{\ast }, p^{\ast }, q^{\ast }$ and $r^{\ast }$).
\end{enumerate}  
\end{description}

\vspace{4mm}\noindent\textbf{Remarks.} We see that, when $v=0$, and $a,u>0$, the right hand side of  (\ref{except}) yields the
(unnormalized) density of the beta law $\beta _{a,u}=BH(a,u,p,0,r)$  where $p$ and $r$ are arbitrary.  By a proper choice of $p$ and $r$ we can ensure that $(a,u,p,0,r)\in P.$ This shows that beta distributions are BH distributions in the sense of Remark 2 in Section 2.2.

 Next notice that if the right hand side of (\ref{except}) is positive in the whole interval $(0,1)$ it has
to be $0<v<1$; if moreover $a,b>0$ then $h$ is integrable and can be
normalized to become the  density of $BH(a,u,-1,rv,r)$ where $r$ is arbitrary. By taking $r$  large enough we can ensure $(a,u,-1,rv,r)\in P.$ The corresponding law will be called
quasibeta and indicated by $q\beta _{a,u,v}$. Thus beta and quasi-beta distributions are beta hypergeometric according to the definition given at the end of Section 2.2.

\vspace{4mm}\noindent \textbf{Proof.} By assumption for $x\in (0,1)$ 
\begin{equation*}
x^{a-1}\left( 1-x\right) ^{b-1}{}_{2}F_{1}(p,q;r;x)=Cx^{a^{\ast }-1}\left(
1-x\right) ^{b^{\ast }-1}{}_{2}F_{1}(p^{\ast },q^{\ast };r^{\ast };x)\text{.}
\end{equation*}%
Dividing the r.h.s. by $x^{a-1}$ and going to the limit as $x\downarrow 0$
we get that the l.h.s. converges to $1$. Then the same has to hold for the
r.h.s., which implies necessarily that $a^{\ast }=a$, and $C=1$. We thus get%
\begin{equation*}
\left( 1-x\right) ^{d}{}_{2}F_{1}(p,q;r;x)=\text{ }_{2}F_{1}(p^{\ast
},q^{\ast };r^{\ast };x)
\end{equation*}%
where $d=b-b^{\ast }$. Recall that $_{2}F_{1}(p,q;r;x)=z(x)$ is a
solution of the second order differential equation (\ref{euler2}).
In the sequel we fix the three numbers $p,q,r$ and for any real number $d$
we define the linear differential operator $L_{d}\left( y\right) \left(
x\right) $ defined by 
\begin{equation}
x\left( 1-x\right) y^{\prime \prime }+[r+\left( 2d-p-q-1\right) x]y^{\prime
}+\frac{dr-pq+[d\left( d-p-q\right) +pq]x}{1-x}y.  \label{ld}
\end{equation}%
We now show that $y(x)=$ $\left( 1-x\right)
^{d}{}_{2}F_{1}(p,q;r;x)$ satisfies $L_{d}\left( y\right) \left( x\right) =0.
$ To see this we write 
\begin{eqnarray*}
z &=&\left( 1-x\right) ^{-d}y,\text{ \ }z^{\prime
}=d(1-x)^{-d-1}y+(1-x)^{-d}y^{\prime } \\
z^{\prime \prime } &=&d(d+1)(1-x)^{-d-2}y+2d\left( 1-x\right)
^{-d-1}y^{\prime }+(1-x)^{-d}y^{\prime \prime }
\end{eqnarray*}%
and we carry the results in (\ref{euler2}). On the other hand $y(x)=\
_{2}F_{1}(p^{\ast },q^{\ast };r^{\ast };x)$, hence 
\begin{equation*}
x\left( 1-x\right) y^{\prime \prime }+[r^{\ast }-\left( p^{\ast }+q^{\ast
}+1\right) x]y^{\prime }-p^{\ast }q^{\ast }y=0  
\end{equation*}%
which subtracted from $L_{d}(y)=0$ yields 
\begin{equation*}
\lbrack r-r^{\ast }+\left( 2d+p^{\ast }+q^{\ast }-p-q\right) x]y^{\prime }+%
\frac{dr+p^{\ast }q^{\ast }-pq+[d\left( d-p-q\right) +pq-p^{\ast }q^{\ast }]x%
}{1-x}y=0
\end{equation*}%
which by setting%
\begin{eqnarray*}
A &=&d\left( d-p-q\right) +pq-p^{\ast }q^{\ast },\text{ \ }B=dr+p^{\ast
}q^{\ast }-pq, \\
C &=&2d+p^{\ast }+q^{\ast }-p-q,\text{ \ }D=r-r^{\ast }
\end{eqnarray*}%
is rewritten as 
\begin{equation*}
\left( Cx+D\right) y^{\prime }+\frac{Ax+B}{\left( 1-x\right) }y=0
\end{equation*}%
Let us now suppose that $C=D=0$ but $A$ and $B$ are not both zero. Then $%
y\equiv 0$, which is impossible. Now assume $C=0$ but $D\neq 0$. Thus%
\begin{equation*}
y^{\prime }=-\frac{Ax+B}{D\left( 1-x\right) }y
\end{equation*}%
and thus 
\begin{equation*}
y=\left( 1-x\right) _{2}^{d}F_{1}(p,q;r;x)=e^{\frac{A}{D}x}\left( 1-x\right)
^{\frac{B-A}{D}}.
\end{equation*}%
Define $d_{1}=d-\frac{B-A}{D}$ and $y_{1}=e^{\frac{A}{D}x}$. Then $y_{1}$
satisfies $L_{d_{1}}\left( y_{1}\right) =0$ (using notation of (\ref{ld})).
Since $L_{d_{1}}\left( y_{1}\right) (x)=\frac{P(x)}{1-x}y_{1}(x)$, where $P$
is a polynomial of degree $3$ with leading coefficient $\left( \frac{A}{D}%
\right) ^{2}$ this implies that $A=0$. Therefore $y=\left( 1-x\right) ^{%
\frac{B-A}{D}}$ and we are in the second case. Next let $C\neq 0$ and $D=0$,
in which case%
\begin{equation*}
\frac{y^{\prime }}{y}=-\frac{Ax+B}{Cx\left( 1-x\right) }=\frac{C_{1}}{x}-%
\frac{C_{2}}{1-x}
\end{equation*}%
and therefore $y=C_{3}x^{C_{1}}\left( 1-x\right) ^{C_{2}}$. Since $y\left(
0\right) =1,$ this implies that $C_{1}=0$ so we are in the second case again.
If $C\neq 0$ and $D\neq 0$, then we can write 
\begin{equation*}
y^{\prime }=\frac{\alpha x+\beta }{\left( 1+cx\right) \left( 1-x\right) }y,
\end{equation*}%
where $c=\frac{C}{D}\neq 0$. We now distinguish between the cases $c\neq -1$
and $c=-1$. If $c\neq -1$ then%
\begin{equation*}
\frac{y^{\prime }}{y}=-\frac{cC_{1}}{1+cx}-\frac{C_{2}}{1-x}
\end{equation*}%
from which $y=\left( 1+cx\right) ^{-C_{1}}\left( 1-x\right) ^{C_{2}}$.
Define $d_{2}=d-C_{2}$ and $y_{2}=\left( 1+cx\right) ^{-C_{1}}$. Then $y_{2}$
satisfies $L_{d_{2}}\left( y_{2}\right) =0$. This implies that $C_{1}\left(
C_{1}+1\right) =0$, therefore either $C_{1}=0$ (thus we are in case b)) or $%
C_{1}=-1$ (thus we are in case c)). Next suppose that $c=-1$. Then%
\begin{equation*}
\frac{y^{\prime }}{y}=\frac{\alpha x+\beta }{\left( 1-x\right) ^{2}}=-\frac{%
C_{1}}{1-x}-\frac{C_{2}}{\left( 1-x\right) ^{2}}
\end{equation*}%
from which $y=\left( 1-x\right) ^{C_{1}}\exp \left\{ \frac{C_{2}}{1-x}%
\right\} $. Define $d_{3}=d-C_{1}$ and $y_{3}=\exp \left\{ \frac{C_{2}}{1-x}%
\right\} $. Then $y_{3}$ satisfies $L_{d_{3}}\left( y_{3}\right) =0$ (using
notation of (\ref{ld})). Since $L_{d_{3}}\left( y_{3}\right) (x)=\frac{P(x)}{%
\left( 1-x\right) ^{3}}y_{3}\left( x\right) $, where $P$ is a polynomial
such that $P(1)=C_{2}^{2}$, which implies that $C_{2}=0$, which falls into
case b). Finally suppose that $A=B=C=D=0$. From $D=0$ we get $r^{\ast }=r$.
Furthermore 
\begin{equation*}
A=B\Longleftrightarrow d(r+d-p-q)=0.
\end{equation*}%
Thus either $d=0$: as a consequence $b=b^{\ast }$, moreover $p+q=p^{\ast
}+q^{\ast }$ (from $C=0$) and $pq=p^{\ast }q^{\ast }$ (from $A=0$), which
means $\left\{ p,q\right\} =\left\{ p^{\ast },q^{\ast }\right\} .$ Or $%
d=p+q-r$: as a consequence $b^{\ast }=b+r-p-q$, moreover $2r-p-q=p^{\ast
}+q^{\ast }$ (from $C=0$) and $\left( r-p\right) \left( r-q\right) =p^{\ast
}q^{\ast }$ (from $A=0$), which means $\left\{ r-p,r-q\right\} =\left\{
p^{\ast },q^{\ast }\right\} $. Thus the condition $A=B=C=D=0$ yields case a)
of the Theorem. To complete part b) if $f_{v}(x)=x^{a-1}\left( 1-x\right)
^{c-1}$ we have 
\begin{equation*}
\left( 1-x\right)^{b-c}\  _{2}F_{1}(p,q;r;x)\equiv 1.
\end{equation*}%
Thus define $d_{4}=b-c$ and $y_{4}\equiv 1$. Then $y_{4}$ satisfies $%
L_{d_{4}}\left( y_{4}\right) =0$ , which implies $d_{4}r-pq=0$ and $%
d_{4}\left( d_{4}-p-q\right) +pq=0$. Summing the two equalities we get $%
d_{4}\left( d_{4}+r-p-q\right) =0$. Thus either $d_{4}=0$ (and $c=b$) and $%
pq=0$; or $d_{4}+r-p-q=0$ (and $c=b+r-p-q$) and $\left( r-p\right) \left(
r-q\right) =0$. Finally, to complete part c), if $f_{v}\left( x\right)
=x^{a-1}\left( 1-x\right) ^{b_{1}-1}\left( 1+cx\right) $, then 
\begin{equation*}
(1-x)^{b_{1}}\left( 1+cx\right) =\left( 1-x\right)^{b}\  _{2}F_{1}\left(
p,q;r;x\right) .
\end{equation*}%
Define $d_{5}=b-b_{1}$ and $y_{5}\left( x\right) \equiv 1+cx$. Then $y_{5}$
satisfes $L_{d_{5}}\left( y_{5}\right) =0$. This means%
\begin{equation*}
c(r+\left( 2d_{5}-p-q-1\right) x)+\frac{d_{5}r-pq+(d_{5}\left(
d_{5}-p-q\right) +pq)x}{1-x}\left( 1+cx\right) =0.
\end{equation*}%
Since $c\neq -1$,$0,$ this implies $d_{5}\left( d_{5}+r-p-q\right) =0.$ As a
consequence the fractional term is the constant $d_{5}r-pq$, so by equating
the coefficients of the polynomial at the l.h.s. to zero we get the two
equations%
\begin{equation*}
\left\{ 
\begin{array}{c}
c\left( 2+r\right) d_{5}=\left( p+1\right) \left( q+1\right)  \\ 
\left( c+d_{5}\right) r=pq%
\end{array}%
\right. 
\end{equation*}%
Thus either $d_{5}=0$, in which case $b_{1}=b$, and $\left( p+1\right)
\left( q+1\right) =0$ and $cr=pq$. Thus if, say, $p=-1$, we get $q=-cr$ as
stated in the theorem. Or $d_{5}=p+q-r$, in which case $b_{1}=b+r-p-q$, and
we get similarly $\left( r-p+1\right) \left( r-q+1\right) =0$ and $cr=\left(
r-p\right) \left( r-q\right) $. Thus if, say, $r-p=-1$, we get $r-q=-cr$ as
stated in the theorem. $\square$
\subsection{Examples} To illustrate the above results, let us recall a few classical equalities (see Abramovitz and Stegun (1965) pages 556-7)
\begin{equation*}\label{HA}
\ _2F_1(p,-p;\frac{1}{2};\sin^2\theta)=\cos 2p\theta,\ \ _2F_1(p,1-p;\frac{1}{2};\sin^2\theta)=\frac{\cos (2p-1)\theta}{\cos \theta}\end{equation*}

\begin{equation*}\label{HB}\ _2F_1(p,1-p;\frac{3}{2};\sin^2\theta)=\frac{\sin (2p-1)\theta}{(2p-1)\sin \theta},\ 
 \ _2F_1(1+p,1-p;\frac{3}{2};\sin^2\theta)=\frac{\sin 2p\theta}{p\sin 2\theta}
\end{equation*}

\begin{equation*}\label{HC}\ _2F_1(p,\frac{1}{2}+p;1+2p;\sin^2\theta)=\frac{1}{(\cos \frac{\theta}{2})^{2p}},\ 
\ _2F_1(p,\frac{1}{2}+p;2p;\sin^2\theta)=\frac{1}{\cos\theta(\cos \frac{\theta}{2})^{2p-1}}
 \end{equation*}
\begin{equation*}\label{HD}\ _2F_1(\frac{1}{2},\frac{1}{2};\frac{3}{2};\sin^2\theta)=\frac{\theta}{\sin \theta},\ 
 \ _2F_1(1,1;\frac{3}{2};\sin^2\theta)=\frac{ 2\theta}{\sin 2\theta}
\end{equation*}
\begin{equation*}\label{HE}\ _2F_1(p,\frac{1}{2}+p;\frac{1}{2};x)=\frac{1}{2(1+\sqrt{x})^{2p}}+\frac{1}{2(1-\sqrt{x})^{2p}},\end{equation*}

\begin{equation*}\label{HF}\ _2F_1(p,\frac{1}{2}+p;\frac{3}{2};x)=\frac{1}{2(1-2p)\sqrt{x}}(\frac{1}{(1-\sqrt{x})^{1-2p}}-\frac{1}{(1+\sqrt{x})^{1-2p}}) 
\end{equation*}
\begin{equation*}\label{HG}\ _2F_1(1,\frac{1}{2};\frac{3}{2};x)=\frac{1}{2\sqrt{x}}\log\frac{1+\sqrt{x}}{1-\sqrt{x}},\ \ _2F_1(1,1;2;x)=\frac{1}{x}\log\frac{1}{1-x},\end{equation*}

Many of these identities  describe a hypergeometric function in terms of $x=\sin^2\theta$: this is sometimes useful for describing a distribution on $(0,\pi/2)$: the image of $BH(a,b,p,q,r)$
by $x\mapsto \theta=\arcsin \sqrt{x}$ from $(0,1)$ to $(0,\pi/2)$ has a density proportional to $\sin^{2a-1}\theta\cos^{2b-1}\theta \ _2F_1(p,q;r;\sin^2 \theta).$

\section{Proof of the basic identity. The $\theta$ parameterization }
\textsc{Mellin-like transforms of $BH(v).$} Assume now that $v=(a,b,p,q,r)^T$ is in $P.$  Recall that this set $P$ has been defined at the very end on Section 2.2 and is 
\begin{equation}\label{PPPP}P=\{(a,b,p,q,r);\ (p,q,r)\in \mathbb{P},\ r\not\in -\mathbb{N}, a,b,r+b-p-q>0\}.\end{equation}
(the set $\mathbb{P}$ is described in the statement of Theorem 2.2).  Let us fix $(s,t)$ such that  $(a+s,b+t,p,q,r)$ is in $P$ or equivalently such that $s>-a$ and $$t>-b-\min \{0,r-p-q\}.$$ 
We get from  (\ref{term}) the important formula 
\begin{equation}
\int_{0}^{1}x ^{s}\left( 1-x\right) ^{t}BH(v)(dx)=
\frac{B(a+s,b+t)\ _{3}F_{2}\left( p,q,a+s;r,a+b+s+t;1\right) }{
B(a,b)\ _{3}F_{2}\left( p,q,a;r,a+b;1\right) }  \label{arimellin}
\end{equation}%
\textsc{The Thomae formula.} There is a fundamental relation between the Euler's gamma function $\Gamma $
and the generalized hypergeometric function $_{3}F_{2}$ defined by (\ref{TFD}) and  evaluated at $x
=1$. This relation has been originally obtained by\ Thomae (1879). In the paper by Maier (2005) and
Beyer \textit{et al }(1987)  such a relation is reformulated as the
invariance of a suitably defined function with respect to the symmetric group $%
\mathcal{S}_{5}$. 

\vspace{4mm}\noindent\textbf{Lemma 3.1.} (Thomae's formula) The function 
\begin{equation*}T(A,B,C,D,E)=
\frac{\Gamma \left( C\right) \  _{3}F_{2}\left(
A,B,C;D,E;1\right) }{\Gamma \left( D\right) \Gamma \left( E\right) }
\end{equation*}%
has the invariance property $T(A,B,C,D,E)=$
\begin{equation*}
T\left( D-C,E-C,D+E-A-B-C,D+E-A-C,D+E-B-C\right)
\end{equation*}

\vspace{4mm}\noindent \textbf{Proof.} See Bailey (1935), Askey, Andrews and Roy (2000) for two
different proofs and Asci \textit{et al} (2008) for a probabilistic one based on the following idea: if $U$ and $V$ are two arbitrary beta random variables on $(0,1)$ and if $t$ is real, compute $\mathbb{E}((1-UV)^t)$ in two ways: expansion in a series of powers of $UV$ or computation of the density of $UV$ by multiplicative convolution. We obtain in this way an identity involving the parameter $t$ and the 4 parameters of the beta distributions. This identity is equivalent to the Thomae's formula. $\square$

\vspace{4mm}\noindent
It should be emphasized that Thomae's formula is an equality between
analytic functions in their whole domain of analyticity, so certainly holds
true when all the arguments $A,B,C,D-C,E-C$ and $D+E-A-B-C$ are positive.

\vspace{4mm}\noindent \textsc{Proof of the basic identity (\ref{BI}).}

\vspace{4mm}\noindent \textbf{Theorem 3.2} Let $v=(a,b,p,q,r)^T$ be in $P$ and assume $a+b-p>0$ and $r-a>0.$ Consider
\begin{equation}\label{TM}M=\left[\begin{array}{rrrrr}-1&0&0&0&1\\1&1&-1&0&0\\0&1&-1&-1&1\\0&1&0&0&0\\0&1&0&-1&1\end{array}\right],\ \Pi=\left[\begin{array}{rrrrr}1&1&-1&0&0\\-1&0&0&0&1\end{array}\right]\end{equation}
Then $Mv$ is in $P$ and if $X\sim BH(v)$ and $W\sim \beta^{(2)}_{\Pi v}$ we have $\frac{1}{1+XW}\sim BH(Mv).$

\vspace{4mm}\noindent \textbf{Proof:} Denote $Mv=(a',b',p',q',r')^T.$ The facts that $v$ is in $P,$ $b'=a+b-p>0$ and $a'=r-a>0$ imply that $p'=r+b-p-q>0$ and
$q'=b>0.$ It is more delicate to prove that $r'=r+b-q>0.$ This certainly holds
when $q \leq 0$ and, from $r+b-p-q>0$, it holds when $p\geq 0$ as
well. Next we assume both $p<0$ and $q>0$. Since $(p,q,r)\in
\mathbb{P}$ we make the following deductions. Case 1 is impossible.
In Case 2 from $r-q \geq 0$ we get $r+b-q>0$. In Case 3 we have
$(p,q,r)\in S$ and since $q>0$ and $r>a>0$ it has to be $1/\Gamma(p)
\geq 0$ which is fulfilled with $p<0$ only when $p=-n$, n being a
positive integer, which is excluded since it would imply $r-p>n$, a
contradiction with $r-p\leq 1$. Finally, in Case 4 we have $0<r<1,
0<q \leq 1$ and $1/\gamma(r-q) \geq 0$. Thus either $r-q>0$, which
implies $r+b-q>0$, or $-1<r-q\leq 0$. The only possibility is thus
$r=q$, in which $r+b-q>0$ as well. ).
Finally, the inequality $r'+b'-p'-q'=a>0$
together with the previous ones implies $Mv \in P$ (see the remark
at the end of Section 2.2).

 Now denote $Y=\frac{1}{1+XW}.$ This implies 
$\frac{1-Y}{Y}=XW$ and therefore for $t\in (-b,a)$ we can write $\mathbb{E}((1-Y)^tY^{-t})=\mathbb{E}(X^t)\mathbb{E}(W^t).$ Because in the following calculation the constants are quite long to write, let us adopt the following convention: we say that two positive functions $f$ and $g$ of $t\in (-b,a)$ are equivalent if $t\mapsto f(t)/g(t)$ is a constant with respect to $t.$ This fact is denoted $f\equiv g$ or -with some abuse of notation- $f(t)\equiv g(t).$ With this convention we get by replacing $(s,t)$ in (\ref{arimellin}) by $(t,0)$
$$\mathbb{E}(X^t)\equiv \frac{\Gamma(a+t)}{\Gamma(a+b+t)}\ _{3}F_{2}( p,q,a+t;r,a+b+t;1)$$
as well as $$\mathbb{E}(W^t)\equiv \Gamma(a+b-p+t)\Gamma(r-a-t).$$ This implies 

\begin{eqnarray} \mathbb{E}((1-Y)^tY^{-t})&\equiv&  \frac{\Gamma(a+b-p+t)\Gamma(r-a-t)\Gamma(a+t)}{\Gamma(a+b+t)}\label{TSP}\\&&\times\ _{3}F_{2}( p,q,a+t;r,a+b+t;1)\nonumber\end{eqnarray}
From formula (\ref{arimellin}) note  that 
\begin{eqnarray}\nonumber \int_{0}^{1}y ^{t}\left( 1-y\right) ^{-t}BH(Mv)(dx)&\equiv&\Gamma(a'+t)\Gamma(b'-t) \\&&\times \ _{3}F_{2}( p',q',a'+t;r',a'+b';1)\nonumber\\&\equiv&\Gamma(r-a+t)\Gamma(a+b-p-t)\label{TSS} \\&&\times \ _{3}F_{2}( r+b-p-q,b,r-a+t;r+b-q,r+b-p;1)\nonumber\end{eqnarray}
The knowledge of the function $t\mapsto \mathbb{E}((1-Y)^tY^{-t})$ gives the knowledge of the distribution of $(1-Y)/Y$ and  of the distribution of $Y.$ Therefore enough is to show that the right hand sides of  $\ref{TSP}$ and $\ref{TSS}$
are equivalent. To see this we simply apply the Thomae formulae (Lemma 3.1) to
$$A=p,\ B=q,\ C=a+t,\ D=a+b+t, E=r$$ and this concludes the proof of Theorem 3.2. $\square$

\vspace{4mm}\noindent
\textsc{The $\theta$ parameterization.}
Up to now, the $BH$ distributions have been parameterized by $(a,b,p,q,r)$ belonging to the subset $P$ of $\mathbb{R}^5$ described in (\ref{PPPP}). One defect of this parameterization is the fact that  $$BH(a,b,p,q,r)=B(a,r+b-p-q,r-p,r-q,r)$$  (as implied by the Euler formula) is not apparent. A second defect of the parameterization $(a,b,p,q,r)$ is that it makes complicated the application of the basic identity. For these reasons we choose to make a linear transformation of $a,b,p,q,r$ as follows. 
Introduce a $5$-tuple of parameters $\mathbf{\theta }=\left( \theta
_{1},\theta _{2},\theta _{3},\theta _{4},\theta _{5}\right) \in \mathbb{R}%
^{5}$. For notational convenience, we  set 
\begin{eqnarray}
a&=&\theta_{4}+\theta_5,\ \  b=\theta_{1}+\theta_4,\nonumber \ \   p=\theta_{3}+\theta_5,\ \  q=\theta_{3}+\theta_4,\\  r&=&\theta_{2}+\theta_3+\theta_{4}+\theta_5
\label{change}
\end{eqnarray}
This can be inverted as 
\begin{eqnarray}2\theta_1=a+2b-p-q,& 2\theta_2=-a-p-q+2r,&\nonumber\\ 
2\theta_3=-a+p+q&\ 2\theta_4=a-p+q,&  2\theta_5=a+p-q\label{ITT}\end{eqnarray}
From now on we denote $BH_{\theta}=BH(a,b,p,q,r)$ where $(a,b,p,q,r)$ is given by (\ref{change}).   For example  important particular cases like $BH(1,1,1,1,2)$ considered in (\ref{BHone}) and $BH(2,2,2,2,4)$ are rewritten as 
\begin{equation}\label{BHone}BH_{\frac{1}{2},\frac{1}{2},\frac{1}{2},\frac{1}{2},\frac{1}{2}}(dx)=\frac{6}{\pi^2 x}\log\frac{1}{1-x}\textbf{1}_{(0,1)}(x)dx.\end{equation}
\begin{equation}\label{BHtwo}BH_{1,1,1,1,1}(dx)=\frac{2}{10-\pi^2 }\left[\frac{2-x}{x^3}\log\frac{1}{1-x}-\frac{2}{x^2}\right]\textbf{1}_{(0,1)}(x)dx.\end{equation}
We say that  $BH_{\theta}$ exists when $(a,b,p,q,r)\in P.$ Necessary and sufficient conditions for this are given in Proposition 3.4 below. This new parameterization has many advantages.
We see immediately that exchange of $p$ with $q$ is equivalent to the
exchange of $\theta _{4}$ with $\theta _{5}$, whereas, since%
\begin{equation*}
r+b-p-q=\theta_{1}+\theta_2,\ r-p=\theta_{2}+\theta_4,\ 
r-q=\theta_{2}+\theta_5 
\end{equation*}%
Euler's identity corresponds to exchange $\theta _{2}$ with $\theta _{3}$
and $\theta _{4}$ with $\theta _{5}$. Therefore actually the distribution \begin{equation}\label{SFI}BH_{\theta
_{1},\theta _{2},\theta _{3},\theta _{4},\theta _{5}}=BH_{\theta
_{1},\{\theta _{2},\theta _{3}\},\{\theta _{4},\theta _{5}\}}\end{equation}
being symmetric in $(\theta_2,\theta_3)$ and $(\theta_4,\theta_5)$ has rather to be considered as depending on $\theta_1$ and on the two sets $\{\theta _{2},\theta _{3}\},\{\theta _{4},\theta _{5}\}.$ This notation $BH_{\theta
_{1},\{\theta _{2},\theta _{3}\},\{\theta _{4},\theta _{5}\}}$ is however a slight abuse of language since the set $\{\theta_2,\theta_3\} $  could be reduced at one point if $\theta_2=\theta_3$ and the language of multisets (sets with entire positive weights) could be more adapted. Up to this, we consider that the notations (\ref{SFI}) are sufficiently informative. 
The revisited basic identity  can be reformulated in this new notation as follows: 

\vspace{4mm}\noindent \textbf{Theorem 3.3.} Let $X\sim BH_{\theta
_{1},\{\theta _{2},\theta _{3}\},\{\theta _{4},\theta _{5}\}}$, $W\sim \beta^{(2)}_{\theta_1+\theta_5,\theta_2+\theta_3}$ and $W'\sim \beta^{(2)}_{\theta_1+\theta_4,\theta_2+\theta_3}$ such that $(W,W')$ are independent of $X.$ Then $$\frac{1}{1+WX}\sim  BH_{\theta
_{5},\{\theta _{1},\theta _{4}\},\{\theta _{2},\theta _{3}\}},\ \frac{1}{1+W'X}\sim  BH_{\theta
_{4},\{\theta _{1},\theta _{5}\},\{\theta _{2},\theta _{3}\}}.$$

\vspace{4mm}\noindent
 Theorem 3.3   shows again that there are two ways to apply the basic identity. It shows also that the matrix $M$ appearing in (\ref{TM}) is similar to a permutation matrix of order 5. 

\vspace{4mm}\noindent \textsc{The existence of $BH_{\theta}
.$}  
In order to check whether $BH_{\theta
_{1},\{\theta _{2},\theta _{3}\},\{\theta _{4},\theta _{5}\}}$ does exist we dissymetrize $\{\theta _{2},\theta _{3}\}$ and $\{\theta _{4},\theta _{5}\}$ by assuming $\theta_2\leq \theta_3$ and $\theta _{4}\leq \theta _{5}.$
Recall that the subset $S$ of $\mathbb{R}^3$ has been defined in (\ref{S}) and is the set of $(x,y,z)$ such that $1/\Gamma(x)\Gamma(y)\Gamma(z)\geq 0.$
The  condition $(a,b,p,q,r)\in P$ where $P$ is given by (\ref{PPPP}) gives the following 

\vspace{4mm}\noindent \textbf{Proposition 3.4.} The distribution $BH_{\theta
_{1},\{\theta _{2},\theta _{3}\},\{\theta _{4},\theta _{5}\}},$ where $\theta_2\leq \theta_3$ and $\theta _{4}\leq \theta _{5}$ exists if and only if 
\begin{itemize}
\item $\theta_1+\theta_2>0$ and $\theta_4+\theta_5>0$ and 
\item either $r=\theta_2+\theta_3+\theta_4+\theta_5\geq 1$ and $\theta_3+\theta_4>0$
\item or $r=\theta_2+\theta_3+\theta_4+\theta_5< 1$ , $\theta_2+\theta_5\leq 1$ and $(\theta_3+\theta_4,\theta_3+\theta_5,r)\in S.$
\end{itemize}

\vspace{4mm}\noindent \textbf{Definition.} We will call $\Theta$ the
set of parameters ${\theta _{1},\theta _{2},\theta
_{3},\theta _{4},\theta _{5}} \in \mathbb{R}^5$ which is the image of $P$ by the linear map $(a,b,p,q,r)\mapsto (\theta_1,\ldots,\theta_5)$ described by (\ref{ITT}). The part of the set $\Theta$ such that $\theta_2\leq \theta_3$ and $\theta _{4}\leq \theta _{5}$ is also described by Proposition 3.4.

 \vspace{4mm}\noindent \textbf{Remark.} From
(\ref{change}) one can observe that $a,b,p,q,r,r+b-p-q>0$ is equivalent to 
\begin{equation}\label{fedup}
\theta_1 +\theta_2, \theta_1+\theta_4, \theta_3 +\theta_4, \theta_4
+\theta_5, \theta_2+\theta_3+\theta_4+\theta_5 >0.
\end{equation}%
 Thus the vectors
$\theta$ satisfying the inequalities (\ref{fedup}) belong
necessarily to $\Theta$. This certainly happens if $\theta$ is such
that $\theta_i +\theta_j>0$ for all $1\leq i<j\leq 5$ except for
$(i,j)=(2,4)$ (here $\theta_2\leq \theta_3$ and $\theta _{4}\leq
\theta _{5}$). Thus in this case $BH_{\theta}$ exists (this
observation will turn out to be useful for Theorem 5.1 below).
Moreover, recalling the first part of the proof of Theorem 3.2 we
see that the application of the basic identity always yields beta
hypergeometric distributions with the vector of parameters $\theta$
satisfying the inequalities (\ref{fedup}).

\section{Random continued fractions with a beta hypergeometric distribution}

It is clear the the iteration of the random transformations
appearing in Theorem 3.3, applied to $X\sim BH_{\theta} $ with $\theta =(\theta _{1}, \theta _{2},\theta
_{3},\theta _{4},\theta _{5})\in
\Theta$ yield random variables whose distribution is of the form
$BH_{\theta'}$ where
$\theta' =(\theta' _{1}, \theta' _{2},\theta'_3,\theta' _{4},\theta' _{5})\in \Theta$ is obtained from a permutation of the
components of $\theta$.

For this reason, for any $\theta \in \mathbb{R}^5$ we define
the finite subset $V_{\theta} \subset \mathbb{R}^5_s$ of vectors $\theta^{\ast}$
which can be obtained in this way. Motivated by Theorem 3.3, we are
going to define a directed graph structure on $V_{\theta}$. The
possible forms of these graphs will be quite limited.

\subsection{The graphs $G_{\theta}$ and their subgraphs}
\textsc{The role of the seven partitions of 5}. There are seven partitions of $5$ enumerated in (\ref{PART}). 
 To each point $\theta\in \mathbb{R}^5$ we attach the discrete measure on $\mathbb{R}$ which is $\sum_{j=1}^5\delta_{\theta_j}=\sum_{k=1}^n m_k\delta_{x_k}$  where $\{x_1,\ldots,x_n\}$ is the image of $j\mapsto \theta_j$ and where $m_k$ is the positive integer which is the number of $j=1,\ldots,5$ such that $\theta_j=x_k.$ Thus $n\leq 5$ and $m_1+\ldots+m_n=5$ defines the  partition of $5$ induced by $\theta\in \mathbb{R}^5.$ For convenience in the sequel we take $m_1\geq m_2\geq \ldots \geq m_n$ and we write 
$$x=x_1, \ y=x_2, \ z=x_3, \ u=x_4, \ v=x_5$$
when these $x_k$ do exist.  Suppose for instance that the partition attached to $\theta$ is 3+2. Therefore we shall use three times the letter $x$ and two times the letter $y;$ 
 the 5 elements of $V_{\theta}$ will be
$$(x,\{x,x\},\{y,y\}),\ (y,\{x,y\},\{x,x\}),\ (x,\{y,y\},\{x,x\}),$$$$\ (y,\{x,x\},\{x,y\}),\ (x,\{x,y\},\{x,y\})$$ 
that we quickly code
\begin{equation}\label{XXXYY}x|x^2|y^2,\ y|xy|x^2,\ x|y^2|x^2,\ y|x^2|xy,\ x|xy|xy \end{equation}

\vspace{4mm}\noindent\textsc{The directed graph $G_{\theta}.$} According to Theorem 3.3, given $\theta\in \mathbb{R}^5$, if $\theta\in \Theta$ (which means that $BH_{\theta}$ exists) and if $\theta_1+\theta_5,$ $\theta_2+\theta_3$ and $\theta_1+\theta_4$ are positive we can move to two new beta hypergeometric distributions. Ignore for a while the constraints linked to inequalities. Let us extend the process to the whole $\theta$'s in $\mathbb{R}^5$ or rather to the elements of the quotient described by $(\theta
_{1},\{\theta _{2},\theta _{3}\},\{\theta _{4},\theta _{5}\}):$ from  this element we move either to $(\theta
_{5},\{\theta _{1},\theta _{4}\},\{\theta _{2},\theta _{3}\})$ or to $(\theta
_{4},\{\theta _{1},\theta _{5}\},\{\theta _{2},\theta _{3}\}).$ These two elements may not be distinct. 

We need to introduce in (\ref{XXXYY}) 
the arrows 
\begin{eqnarray*}
x|x^2|y^2&\rightarrow& y|xy|x^2\\
y|xy|x^2&\rightarrow& x|xy|xy \\
x|y^2|x^2&\rightarrow&x|x^2|y^2\\
y|x^2|xy&\rightarrow&x|y^2|x^2,\ y|xy|x^2\\
x|xy|xy&\rightarrow&y|x^2|xy,\ x|xy|xy
\end{eqnarray*}
getting the following directed graph, called the graph $G_{\theta}$: 
\begin{equation}\label{HEX}\xymatrix{ {x|x^2|y^2} \ar@{>}[r]& {y|xy|x^2} \ar@{>}[dr] &\\
 {x|y^2|x^2}\ar@{>}[u]& {y|x^2|xy} \ar@{>}[l]\ar@{>}[u]&  {x|xy|xy}\ar@{>}[l] \ar@(ur,ul)[]}\end{equation}

It is clear that the graph $G_{\theta}$ is the same for all $\theta
\in \mathbb{R}^5$ sharing the same partition of $5$ given by the
numbers $m_1 \geq m_2 \geq \ldots m_n$. Thus we explore the 7
possible forms of this graph. We are particularly interested in
determining the cycles in these graphs. A cycle of order $k\geq 2$
in a directed graph is a sequence $v_0,\ldots,v_{k-1}$ of
\textit{distinct }vertices such that $(v_{k-1},v_0)$ and
$(v_i,v_{i+1})$ are oriented edges of the graph for
$i=0,\ldots,k-2.$ A cycle of order 1 is a vertex $v$ such that
$(v,v)$ is an edge).

\vspace{4mm}\noindent\textsc{Description of the seven  graphs $G_{\theta}$} 
\begin{itemize}\item {\textsc{The case 5.}} Here $\theta=(x,x,x,x,x)$ and  the graph $G_{\theta}$ is rather trivial with one point  and one cycle of order 1

\begin{equation}\label{KEX}\xymatrix{{x|x^2|x^2}\ar@(ur,ul)[]}\end{equation}
\item \textsc{The case 4+1.} Here the graph has three vertices and is 
\begin{equation}\label{JEX}\xymatrix{ {x|x^2|xy} \ar@{>}[r]\ar@{>}[d]& {x|xy|x^2}  \ar@{>}[l]\\
 {y|x^2|x^2}\ar@{>}[ur]& }\end{equation} It has two  cycles of order 2 and 3.
\item \textsc{The case 3+2.} The graph has already been drawn in (\ref{HEX}). It has three  cycles of order 1, 3 and 5.

\item \textsc{The case 3+1+1.} Here the graph has eight vertices and is

\[\xymatrix@-3ex{ &{z|x^2|xy} \ar@{>}[r]\ar@{>}[d]& {x|yz|x^2} \ar@{>}[d] &{y|x^2|xz} \ar@{>}[d]\ar@{>}[l]&\\
 & {y|xz|x^2} \ar@{>}[dl]&  {x|x^2|yz}\ar@{>}[l]\ar@{>}[r]&{z|xy|x^2}\ar@{>}[dr]&\\
 {x|xz|xy}\ar@{>}[uur]\ar@{>}[rrrr]&&&&{x|xz|xy}\ar@{>}[uul]\ar@{>}[llll]
 }\]
There are one  cycle of order 2,  two of orders 3, 5 and 6.

\item \textsc{The case 2+2+1.} Here the graph has 11 vertices and is 
\begin{equation}\label{GEX}\xymatrix@-5ex{ &&&{z|y^2|x^2}\ar@{>}[dlll]&&&\\
             {x|xz|y^2}\ar@{>}[rr]&&{y|xy|xz}\ar@{>}[rr]\ar@{>}[dr]&&{x|yz|xy}\ar@{>}[rr]\ar@{>}[dd]&&{y|x^2|yz}\ar@{>}[dd]\ar@{>}[ulll]\\
             &&&{z|xy|xy}\ar@{>}[dl]\ar@{>}[ur]&&&\\
             {x|y^2|xz}\ar@{>}[uu]\ar@{>}[drrr]&&{y|xz|xy}\ar@{>}[ll]\ar@{>}[uu]&&{x|xy|yz}\ar@{>}[ll]\ar@{>}[ul]&&{y|yz|x^2}\ar@{>}[ll]\\
             &&&{z|x^2|y^2}\ar@{>}[urrr]&&&}\end{equation}
There are two  cycles of order 3,  one of order 4, six of order 5,  four of order 6, two of order 7 and 9 and one of order 8.

\item \textsc{The case 2+1+1+1.} The 18 vertices  are
$$\begin{array}{ccc}A=x|xu|yz &A'=x|xy|uz &A''=x|xz|uy\\
B=x|yz|xu &B'=x|uz|xy &B''=x|uy|xz\\
C=u|x^2|yz &C'=y|x^2|yz &C''=z|x^2|uy\\
D=u|yz|x^2 &D'=y|uz|x^2 &D''=z|uy|x^2\\
E_1=u|xz|xy &E'_1=y|ux|xz &E''_1=z|xy|ux\\
E_2=u|xy|xz &E'_2=y|xz|ux &E''_2=z|ux|xy\\
\end{array}$$
Here is the graph
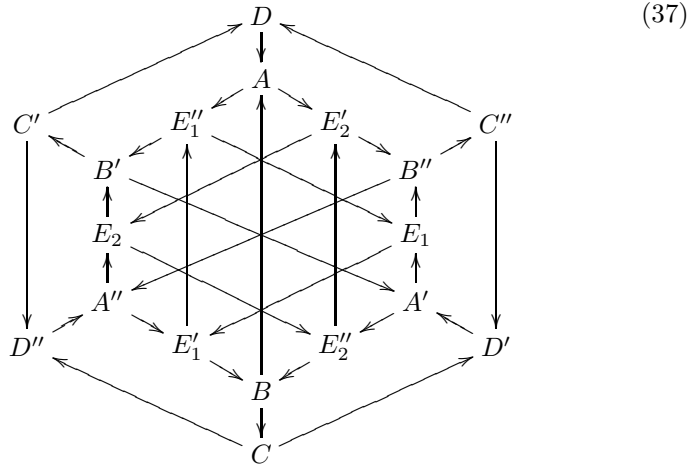
\begin{equation}\label{MEX}\xymatrix@-5ex{ 
&&&&&&D\ar@{>}[dd]&&&&&&\\
&&&&&&&&&&&&\\
&&&&&&A\ar@{>}[dll]\ar@{>}[drr]&&&&&&\\
C'\ar@{>}[uuurrrrrr]\ar@{>}[dddddd]&&&&E''_1\ar@{>}[dll]\ar@{>}[dddrrrrrr]&&&&E'_2\ar@{>}[drr]\ar@{>}[dddllllll]&&&&C''\ar@{>}[uuullllll]\ar@{>}[dddddd]\\
&&B'\ar@{>}[ull]\ar@{>}[ddddrrrrrrrr]&&&&&&&&B''\ar@{>}[urr]\ar@{>}[ddddllllllll]&&\\
&&&&&&&&&&&&\\
&&E_2\ar@{>}[uu]\ar@{>}[dddrrrrrr]&&&&&&&&E_1\ar@{>}[uu]\ar@{>}[dddllllll]&&\\
&&&&&&&&&&&&\\
&&A''\ar@{>}[drr]\ar@{>}[uu]&&&&&&&&A'\ar@{>}[dll]\ar@{>}[uu]&&\\
D''\ar@{>}[urr]&&&&E'_1\ar@{>}[drr]\ar@{>}[uuuuuu]&&&&E''_2\ar@{>}[dll]\ar@{>}[uuuuuu]&&&&D'\ar@{>}[ull]\\
&&&&&&B\ar@{>}[dd]\ar@{>}[uuuuuuuu]&&&&&&\\
&&&&&&&&&&&&\\
&&&&&&C\ar@{>}[uuurrrrrr]\ar@{>}[uuullllll]&&&&&&}\end{equation}
There are two  cycles of order 3,  twelve of order 5, nine of order 6,  three of order 7, nine of   order  8,  eight of order 9, three of order 10, three  of order 12,  six  of order 13, six   of order  14,  two of order 15, nine of   order  16.

\item \textsc{The case 1+1+1+1+1.} The graph has 30 vertices and is too complicated to be drawn here. The two edges issued from  $u|vx|yz$ are given by 
$u|vx|yz\rightarrow y|uz|vx, z|uy|vx.$ There are exactly two incoming edges, coming  from $x|yz|uv$ and $v|yz|ux.$ There is a large number of cycles in this graph; the following remark helps in their determination.  
\end{itemize}

\vspace{4mm}\noindent\textsc{A remark about the automorphisms of the graphs and their cycles.}
The graphs 3+1+1, 2+2+1, 2+1+1+1+1 and 1+1+1+1+1 have automorphisms induced by the permutations of the letters. For instance, the vertices of the graph 2+1+1+1 are coded by letters $x^2yzu$ and the 6 permutations of $yzu$ induce a group $G$ of automorphisms of the graph. Clearly, $G$ transforms a  cycle of size $k$ into a  cycle of size $k.$ Therefore the set of  cycles of size $k$ is splitted in orbits. For the simpler graphs $3+1+1,\ 2+2+1,$ and  $2+1+1+1$ the number of orbits  can be easily found by hand. We indicate below the number of orbits of size $k$ for the graph 1+1+1+1+1, which have been determined by computer. We have  not displayed  the sometimes quite large number of  cycles of each order  as we did for the six others. 
We get the following results
\begin{eqnarray*}
3+1+1&:&2(1),\ 3(1),\ 5(1),\ 6(1)\\
2+2+1&:&3(1),\ 4(1),\ 5(3),\ 6(2),\ 7(1),\ 8(1),\ 9(1)\\
2+1+1+1&:&3(1),\ 5(2),\ 6(2),\ 7(1),\ 8(2),\ 9(2),\\&& 10(1),\ 12(1),\ 13(1),\ 14(1),\ 15(1),\ 16(3)\\ 
1+1+1+1+1&:&5(1),\ 6(1),\ 8(1),\ 9(1),\ 12(2),\ 13(1),\\&& 14(3),\ 15(4),\ 16(7),\ 17(3),\ 
18(4),\ 19(8),\\&& 20(7),\ 22(7),\ 23(10),\ 24(2),\ 26(15),\ 30(4)
\end{eqnarray*}
To understand this array, 30(4) on the last line means that the graph 1+1+1+1+1 has 4 different orbits on the set of cycles of order 30 (the existence of cycles of order 30 implies that the graph is Hamiltonian). 

\vspace{4mm}\noindent\textsc{The  two subgraphs $G^*_{\theta}$ and $G^{**}_{\theta}$ of $G_{\theta}$}

Consider such a graph  $G=G_{\theta}$ that we have just defined. 
Denote  \begin{equation}\label{VVV}v_{\theta}=(\theta_1,\{\theta_2,\theta_3\},\{\theta_4,\theta_5\})\end{equation} and write  $v_0=v_{\theta}$. If $(v_0,v)$ is an edge of the graph recall that either 
$v=v_1=(\theta_5,\{\theta_1,\theta_4\},\{\theta_2,\theta_3\})$ or $v=v_2=(\theta_4,\{\theta_1,\theta_5\},\{\theta_2,\theta_3\}).$ We say that  that the edge $ (v_0,v)$ is acceptable if 
$\theta_2+\theta_3>0$ and either $\theta_1+\theta_5>0$ when $v=v_1$ or $\theta_1+\theta_4>0$ when $v=v_2.$ We denote by $G^*_{\theta}$ the subgraph of $G$ when we remove the non acceptable edges. 
Finally, in the graph $G^*_{\theta}$ let us remove the vertices $v$ such that $BH_{v}$ does not exist. We also remove the edges of $G^*_{\theta}$ which are adjacent to these erased vertices. The remaining graph is denoted by $G^{**}_{\theta}.$ A detailed example is in order: we start from $\theta=z|y^2|x^2$ where $x,y,z$ are distinct real numbers. Therefore $G_{\theta}$ is of the 2+2+1 type and it is the graph (\ref{GEX}). We now assume $x+y>0,y+z>0, y<0.$ This leads to the following graph $G^*_{\theta}$
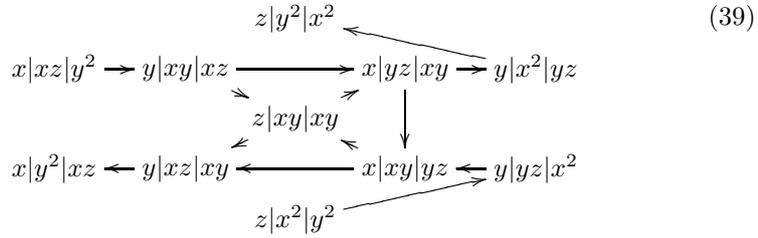
\begin{equation}\label{GEY}\xymatrix@-5ex{ &&&{z|y^2|x^2}&&&\\
         {x|xz|y^2}\ar@{>}[rr]&&{y|xy|xz}\ar@{>}[rr]\ar@{>}[dr]&&{x|yz|xy}\ar@{>}[rr]\ar@{>}[dd]&&{y|x^2|yz}\ar@{>}[ulll]\\
             &&&{z|xy|xy}\ar@{>}[dl]\ar@{>}[ur]&&&\\
             {x|y^2|xz}&&{y|xz|xy}\ar@{>}[ll]&&{x|xy|yz}\ar@{>}[ll]\ar@{>}[ul]&&{y|yz|x^2}\ar@{>}[ll]\\
             &&&{z|x^2|y^2}\ar@{>}[urrr]&&&}\end{equation}
             We assume furthermore $x+\min (x,z)+2y\geq 1$ and we get $G^{**}_{\theta}$  
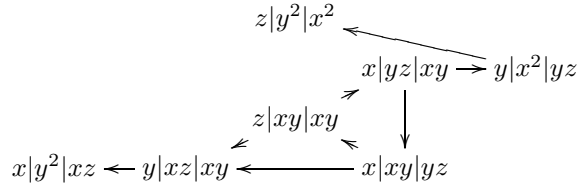
\begin{equation*}\label{GEZ}\xymatrix@-5ex{ &&&{z|y^2|x^2}&&&\\
         &&&&{x|yz|xy}\ar@{>}[rr]\ar@{>}[dd]&&{y|x^2|yz}\ar@{>}[ulll]\\
             &&&{z|xy|xy}\ar@{>}[dl]\ar@{>}[ur]&&&\\
             {x|y^2|xz}&&{y|xz|xy}\ar@{>}[ll]&&{x|xy|yz}\ar@{>}[ll]\ar@{>}[ul]&&
             }\end{equation*}

Here is a second example with $\theta=x|x^2|y^2$ where $x+y\geq 1/2 $ and $y<0$. Here $G_{\theta}$ is of the 3+2 type and the subgraphs  $G^*_{\theta}$ and $G^{**}_{\theta}$  are respectively

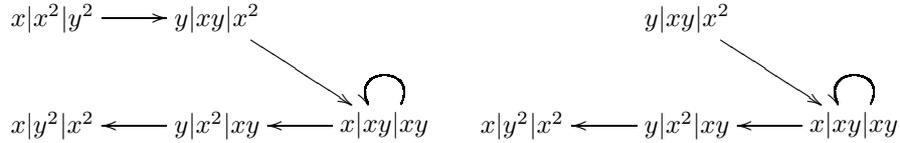
\begin{equation*}\label{HEY}\xymatrix{ {x|x^2|y^2} \ar@{>}[r]& {y|xy|x^2} \ar@{>}[dr] &\\
 {x|y^2|x^2}& {y|x^2|xy} \ar@{>}[l]&  {x|xy|xy}\ar@{>}[l]\ar@(ur,ul)[] }\ \ \ \ \xymatrix{ & {y|xy|x^2} \ar@{>}[dr] &\\
 {x|y^2|x^2}& {y|x^2|xy} \ar@{>}[l]&  {x|xy|xy}\ar@{>}[l] \ar@(ur,ul)[]}\end{equation*}

The  third example with $\theta=z|x^2|y^2$ with $y>x>0,$ $y>1$ and $z=-x$ gives the two graphs $G^{*}_{\theta}$ and $G^{**}_{\theta}$.  The graph $G^{**}_{\theta}$ has no cycle at all. 
$$
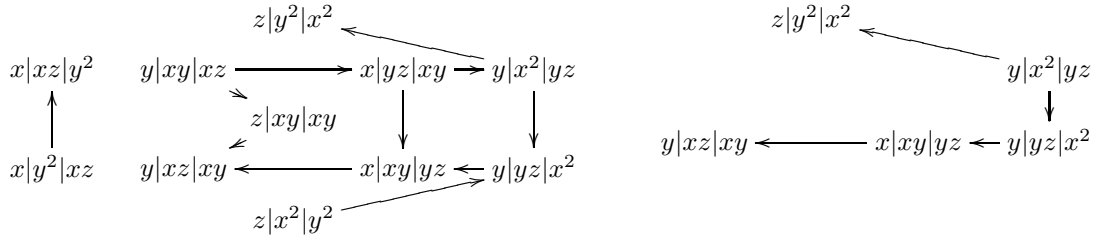
\xymatrix@-5ex{ &&&{z|y^2|x^2}&&&\\
             {x|xz|y^2}&&{y|xy|xz}\ar@{>}[rr]\ar@{>}[dr]&&{x|yz|xy}\ar@{>}[rr]\ar@{>}[dd]&&{y|x^2|yz}\ar@{>}[dd]\ar@{>}[ulll]\\
             &&&{z|xy|xy}\ar@{>}[dl]&&&\\
             {x|y^2|xz}\ar@{>}[uu]&&{y|xz|xy}&&{x|xy|yz}\ar@{>}[ll]&&{y|yz|x^2}\ar@{>}[ll]\\
             &&&{z|x^2|y^2}\ar@{>}[urrr]&&&}
\ \ \ \ \xymatrix@-5ex{ &&&{z|y^2|x^2}&&&\\
             &&&&&&{y|x^2|yz}\ar@{>}[dd]\ar@{>}[ulll]\\
             &&&&&&\\
             &&{y|xz|xy}&&{x|xy|yz}\ar@{>}[ll]&&{y|yz|x^2}\ar@{>}[ll]}$$

\vspace{4mm}\noindent \textbf{Remark.} If $\theta_{i}+\theta_j>0$ for all $1\leq i<j\leq 5$ then $G_{\theta}=G^{**}_{\theta}.$

\subsection{Random continued fractions  attached to a path in $G^{**}_{\theta}$}

\vspace{4mm}\noindent\textsc{The basic identity and the graphs.}
Let us fix $\theta\in \mathbb{R}^5 $ and consider the directed graph $G^{**}_{\theta}=(V,E)$ 
To each vertex $v\in V$ is attached a distribution 
$BH_v$. To each edge $(v,v')\in E$ is attached a pair of positive numbers corresponding to a $\beta^{(2)}$ distribution that we denote by $\beta^{(v,v')}$. 
The basic identity (Theorem 3.3) says that if $X\sim BH_v$ and $W\sim \beta^{(v,v')}$ are independent then $\frac{1}{1+XW}\sim  BH_{v'}.$

In the sequel, for $w>0$, we denote by $H_w$ the  Moebius transformation 
$$H_w(x)=\frac{1}{1+wx}.$$

\vspace{4mm}\noindent \textbf{Proposition 4.1.} Let $v_0\rightarrow v_1\rightarrow\ldots\rightarrow v_n$
be a path in $G^{**}_{\theta}$ of  non necessarily distinct vertices. Let  $X_0,W_1,\ldots,W_n$ be independent random variables such that $X_0\sim BH_{v_0}$ and $W_j\sim \beta^{(v_{j-1},v_{j})}$ for $j=1,\ldots, n$. Define the random Moebius transformations
$F_j=H_{W_j}.$ Then 
$$F_n\circ F_{n-1}\circ\ldots F_1(X_0)=\frac{1}{1+\frac{W_n}{\ldots+\frac{W_2}{1+W_1X_0}}}\sim BH_{v_n}.$$

\vspace{4mm}\noindent \textbf{Proof.} We proceed by induction on $n.$ This is trivially true for $n=0.$ If it is true for $n-1$ we apply the basic identity (Theorem 3.3) to the pair $(v_{n-1},v_n).$$\square$

Here is the simple theorem which can be considered as the main result of the present paper:

\vspace{4mm}\noindent \textbf{Theorem  4.2.} Let 
\begin{equation}\label{IPV}\ldots v_n \rightarrow v_{n-1}\rightarrow\ldots\rightarrow v_2\rightarrow v_1\rightarrow v_0\end{equation}
be an infinite path in $G^{**}_{\theta}$ and let $W_1,\ldots,W_n,\ldots$ be independent random variables such that  $W_n\sim \beta^{(v_{n},v_{n-1})}$. Define the random Moebius transformations
$F_j=H_{W_j}$ and 
$Z_n(x)=F_1\circ F_{2}\circ\ldots F_{n}(x).$
Then the random continued fraction  $$Z=\lim_{n\rightarrow \infty} Z_n(x)=\frac{1}{1+\frac{W_1}{1+\frac{W_2}{1+W_3\ldots}}}$$ associated to the infinite path
almost surely exists  and is independent of $x.$ Its distribution is $Z\sim BH_{v_0}$

\vspace{4mm}\noindent \textbf{Proof.} Since the graph $G^{**}_{\theta}$ is finite, there exists an edge $(v,v')$ such that the set $N=\{n; (v_{n},v_{n-1})=(v,v')\}$ is infinite. Since for all $n\in N$ all the $W_n$ are identically distributed we have $\sum_{n\in N}1/W_n=\infty $ almost surely and therefore $\sum_{n=1}^{\infty}1/W_n=\infty $ almost surely. This is a classical criteria (see Chrystal  (1964)) for claiming that  $Z=\lim_{n\rightarrow \infty} Z_n(x)$ exists and does not depend on $x.$ Let now $X$ be independent of $W_1,\ldots,W_n,\ldots$ such that $X\sim BH_{v}.$ Let us apply Proposition 4.1 to $Z_n(X)$ when $n\in N$: we get that $Z_n(X)\sim BH_{v_0}$ for all $n\in N.$ We deduce from this that $Z\sim BH_{v_0}.$ $\square$

\vspace{4mm}\noindent \textsc{Comments and examples.} In Theorem 4.2 we have called the sequence (\ref{IPV}) an infinite path in the graph $G^{**}_{\theta}$  which can also be written as 
$$v_0\leftarrow v_{1}\leftarrow v_2\leftarrow \ldots \leftarrow v_{n-1}\leftarrow v_n\leftarrow\ldots$$
If one insists that an infinite path should be a map from $\mathbb{N}$ to the set of vertices of the graph and not a map on $-\mathbb{N}$ it would have been be more correct to say that we work with an infinite path in the graph where all arrows have been inverted.

 A simple example of application of Theorem 4.2 is the graph (\ref{KEX}) since obviously the sequence $(v_n)$ is the constant $x|x^2|x^2$ where $x>0.$ Here the distribution of $Z$ is 
$$BH_{x,x,x,x,x}(dz)=Cz^{2x-1}(1-z)^{2x-1}\ _2\ F_1(2x,2x;4x;z)\textbf{1}_{(0,1)}(z)dz$$
Examples for $x=1/2$ and $x=1$ have been given in (\ref{BHone}) and (\ref{BHtwo}).

Another  example of application of Theorem 4.2 is the graph (\ref{HEX}) that we present in a simpler way as 
\begin{equation}\label{JEX}\xymatrix{ {a} \ar@{>}[r]\ar@{>}[d]& {b}  \ar@{>}[l]\\
 {c}\ar@{>}[ur]& }\end{equation}
where $a=x|x^2|xy,\ b= x|xy|x^2,\  c=y|x^2|x^2$. It is easily seen that $x+y, x>0$ implies $G^{**}_{\theta}=G_{\theta}.$  To any infinite word of $\{0,1\}$ we associate an infinite path $ v_0 \leftarrow v_{1}\leftarrow\ldots\leftarrow v_{n-1}\leftarrow v_n\leftarrow \ldots$ in this graph as follows:
each one is replaced by $ b\leftarrow c\leftarrow a\leftarrow$ and each zero is replaced by $b\leftarrow  a\leftarrow.$ For instance the word $00101\ldots$ gives an infinite path  ending at $b$ as
$$b\leftarrow  a\leftarrow b\leftarrow  a\leftarrow b\leftarrow c\leftarrow a\leftarrow b\leftarrow  a\leftarrow   b\leftarrow c\leftarrow a\leftarrow \dots.$$ Theorem 4.5 says that
whatever is the infinite word of $\{0,1\}$ the distribution of the
random continued fraction $Z$ corresponding to the infinite path in
the graph (\ref{JEX}) deduced from this word is
$$BH_b(dz)=Cz^{2x-1}(1-z)^{2x-1}\ _2\
F_1(2x,2x;3x+y;z)\textbf{1}_{(0,1)}(z)dz.$$

 \vspace{4mm}\noindent \textsc{A Cauchy distribution analogy.}
If $w=a+ib$ with $b>0,$  consider the Cauchy distribution $C_w(dx)=\frac{1}{\pi}\frac{bdx}{b^2+(x-a)^2}$. 
 Now let $W_1,\ldots, W_n,\ldots$  be independent random variables such that $W_n\sim C_{w_n}.$ Assume that $\sum_{n=1}^{\infty}|w_n|^{-1}=\infty.$ Define
$$Z=W_1-\frac{1}{W_2-\frac{1}{W_3-\ldots}},\ \ z=w_1-\frac{1}{w_2-\frac{1}{w_3-\ldots}}.$$ Then $Z\sim C_z.$
We find some analogy between this elegant statement (due to Lloyd (1969) in the particular case of a constant sequence $(w_n)_{n= 1}^{\infty}$) and Theorem 4.2: here we consider the Moebius transformations $F_w(x)=
w-\frac{1}{x}$ and we are given an  arbitrary infinite word $w_1w_2\ldots.$ It leads to  the exact distribution $C_z$ of 
$$Z=\lim_{n\rightarrow \infty}F_{W_1}\circ F_{W_2}\circ\ldots F_{W_n}(x).$$

\subsection{Markov chains  attached to a cycle in $G^{**}_{\theta}$}
 A consequence of Theorem 4.2 is about the stationary  distribution of some Markov chains.  It shows the importance of cycles in the graphs $G^{**}_{\theta}.$ Let us  mention that Theorem 5.1 below will prove that the cycles of $G^{*}_{\theta}.$ and $G^{**}_{\theta}$ coincide.

\vspace{4mm}\noindent \textbf{Theorem  4.3.} Let $v_0\rightarrow v_1\rightarrow\ldots v_{k-1}\rightarrow v_0$ be a cycle of order $k$ in $G^{**}_{\theta}$;  let $W_j\sim \beta^{(v_{j-1},v_j)}$ with $j=1,2,\ldots $ be independent with the convention $v_j=v_{j'}$ if $j\equiv j'\ \mathrm{mod}\  k.$ Consider the random Moebius transformation
$$G_n(x)=H_{W_{nk}}\circ H_{W_{nk-1}}\circ \ldots\circ  H_{W_{(n-1)k+1}}(x)$$ and the homogeneous Markov chain $(X_n(x))_{n=0}^{\infty}$ on $(0,1)$ defined by $X_0(x)=x>0$ and $$X_n(x)=G_n(X_{n-1}(x))=G_n\circ G_{n-1}\circ \ldots\circ  G_1(x).$$ Under these circumstances the stationary distribution of the chain is unique and is $BH_{v_0}.$ 

\vspace{4mm}\noindent \textbf{Proof.} Consider $Z_n(x)= G_1\circ G_{2}\circ \ldots\circ  G_n(x).$  Consider the infinite path
$$v_0\leftarrow v_k\leftarrow v_{k-1}\leftarrow\ldots\leftarrow v_1\leftarrow v_0\leftarrow v_k\leftarrow v_{k-1}\leftarrow\ldots
$$ Theorem 4.2 shows that $Z=\lim_{n\rightarrow\infty} Z_n(x)$ exists almost surely and that $Z\sim BH_{v_0}.$ From Letac (1986), from Chamayou and Letac (1991)  or from Propp and Wilson (1995) we get that the stationary distribution of the Markov chain $(X_n(x))_{n=0}^{\infty}$ is unique and it is the distribution of $Z.$   $\square$

\vspace{4mm}\noindent \textsc{Comments and examples.} For other random continued fractions  there are some  analogues of Theorem 4.3 in literature with  cycles only of  size 1 or 2: see   Letac and Seshadri (1984), Lloyd (1959), Dyson (1953),  Marklov, Tourigny and Wolovski (2008) and Asci \textit{et al.} (2008) that we have mentioned in   Section 
4.2.

 The germane example of the present paper is in Asci \textit{et al.} (2008) where the two $BH$ distributions $BH(a,b,a,b,a+a')$ and $BH(a',b,a',b,a+a')$ are considered. Their $\theta$ parameterizations are $$BH_{\frac{b}{2}, \{\frac{b}{2},a'-\frac{b}{2}\},\{\frac{b}{2},a-\frac{b}{2}\}},\ BH_{\frac{b}{2}, \{\frac{b}{2},a-\frac{b}{2}\},\{\frac{b}{2},a'-\frac{b}{2}\}}.$$ These $\theta$'s are of $3+1+1$ type. With $x=b/2$, $y=(2a-b)/2$ and $z=(2a'-b)/2$ the vertices $x|xz|xy$ and $x|xy|xz$ are the vertices of the unique cycle of order 2 of the graph $3+1+1.$ The parameters of the beta type two random variable $W$ used by the random Moebius transformation $H_W$ sending the first law into the second are thus $(b,a)$ (and $(b,a')$ for the opposite).

In Section 4.1 we have mentioned the existence of cycles of all
sizes between 1 and 30 (except 11,\ 21,\ 25,\ 27,\ 28,\ 29). Each of
these cycles  is associated to explicit distributions of periodic random
continued fractions and stationary measures of Markov chains. In
particular we give the example of the homogeneous Markov chain
$G_n\circ G_{n-1}\circ \ldots\circ G_1(z)$ with stationary
distribution $BH_{x,y,z,u,v}$, where the random Moebius
transformations $G_k$ are iid with
$$
G_1=H_{W_{30}}\circ \ldots H_{W_1}
$$
with the $W_i$'s independent with $\beta^{(2)}$ distribution with
parameters (increasing from $i=1$ to $i=30$)
$$
(x+y, u+v), (y+u,x+z), (u+x, v+y), (x+v, z+u), (v+z, x+y),
$$
$$
(z+x, u+v), (x+u, z+y), (u+z, x+v), (z+x, u+y), (x+u, v+z),
$$
$$
(u+v, y+x), (v+y, u+z), (y+u, v+x), (u+v, z+y), (v+z, u+x),
$$
$$
(z+u, v+y), (u+v, z+x), (v+z, y+u), (z+y, x+v), (y+x, z+u),
$$
$$
(x+z, y+v), (z+y, x+u), (y+x, v+z), (x+v, y+u), (v+y, z+x),
$$
$$
(y+z, u+v), (z+u, y+x), (u+y, v+z), (y+v, x+u), (v+x, y+z).
$$

We should also mention here that if $$M=\left[\begin{array}{cc}a&b\\c&d\end{array}\right]$$ is a non singular real matrix and if $h_M(x)=(ax+b)/(cx+d)$ is the corresponding Moebius transformation we have necessarily $h_M\circ h_{M'}=h_{MM'}.$ Therefore the above results could be interpreted in terms of random walks on the group $GL(2,\mathbb{R})$
 suitably quotiented since $h_{\lambda M}=h_M$ for any non zero scalar $\lambda.$

\section{Cycles and positivity} Let $G=G_{\theta}$ be one of the seven graphs.  In this section we show the delicate result that the cycles of $G^*_{\theta}$ are always the cycles of $G^{**}_{\theta}.$ Since the existence of  $\beta^{(v,v')}$ is  easier to check than  the existence of $BH_{v},$ Theorem 5.1 happens to be  useful and practical. 
Use again the notation   \begin{equation*}\label{VVV}v_{\theta}=(\theta_1,\{\theta_2,\theta_3\},\{\theta_4,\theta_5\})\end{equation*} and suppose that  $v_0=v_{\theta}$ is a vertex of $G.$ If $(v_0,v)$ is an edge of the graph recall that either 
$v=v_1=(\theta_5,\{\theta_1,\theta_4\},\{\theta_2,\theta_3\})$ or $v=v_2=(\theta_4,\{\theta_1,\theta_5\},\{\theta_2,\theta_3\})$ and that   $ (v_0,v_1)$ is said to be acceptable if 
$\theta_2+\theta_3>0$ and  $\theta_1+\theta_5>0$.

\vspace{4mm}\noindent\textbf{Theorem 5.1.} Let $\theta=(\theta_1,\theta_2,\theta_3,\theta_4,\theta_5)\in \mathbb{R}^5 $ and consider the vertex $v_0=v_{\theta}$   of $G^*_{\theta}$ defined by (\ref{VVV}). We assume without loss of generality $\theta_2\leq \theta_3$ and $\theta_4\leq \theta_5.$ Then $v_0$ belongs to a cycle of $G_{\theta}^*$ if and only if  $\theta_i+\theta_j>0$ for all $1\leq i<j\leq 5$
except possibly for $(i,j)=(2,4).$ Furthermore if $v_0$ belongs to a cycle of $G_{\theta}^*$ then $v_0\in  \Theta,$ which means that the distribution $BH_{v_0}$ exists and $v_0$ is a vertex of $G^{**}_{\theta}.$
Finally, if $\theta_{2}+\theta_4\leq 0$ then  $$v_0\rightarrow v_1\rightarrow v_2\rightarrow v_3\rightarrow v_4\rightarrow v_5\rightarrow v_0$$where
\begin{eqnarray}\nonumber v_1&=&(\theta_5,\{\theta_4,\theta_1\},\{\theta_2,\theta_3\}),\ v_2=(\theta_3,\{\theta_2,\theta_5\},\{\theta_4,\theta_1\}),\ v_3=(\theta_1,\{\theta_4,\theta_3\},\{\theta_2,\theta_5\}),\\ 
v_4&=&(\theta_5,\{\theta_2,\theta_1\},\{\theta_4,\theta_3\}),\ v_5=(\theta_3,\{\theta_4,\theta_5\},\{\theta_2,\theta_1\}) \label{VCC}\end{eqnarray}
and $\min \{\theta_1,\theta_3,\theta_5\}>-\max\{\theta_2,\theta_4\}.$
In particular in this case $v_0$ belongs to a cycle of order 1,2,3 or 6.

\vspace{4mm}\noindent\textbf{Lemma 5.2.} In the graph $G$, let $v_{-2}\rightarrow v_{-1}\rightarrow v_{0}\rightarrow v_{1}\rightarrow v_{2}$ such that the four edges $\{(v_k,v_{k+1}); \ k=-2,-1,0,1\}$ are acceptable. If $v_0=(\theta_1,\{\theta_2,\theta_3\},\{\theta_4,\theta_5\}) $ then the six following  numbers
$$\theta_1+\theta_2,\ \theta_1+\theta_3,\ \theta_2+\theta_3,\ \theta_1+\theta_4,\ \theta_1+\theta_5,\ \theta_4+\theta_5$$ are positive (we do not assume  $\theta_2\leq \theta_3$ and $\theta_4\leq \theta_5$ here). 

\vspace{4mm}\noindent\textbf{Proof of Lemma 5.2.} The edge $(v_0,v_1)$ being acceptable we get $\theta_2+\theta_3>0.$ The edge $(v_{-1},v_0)$ being acceptable we get $\theta_4+\theta_5>0.$ If $v_1=(\theta_4,\{\theta_1,\theta_5\},\{\theta_2,\theta_3\}) $ the fact that $(v_0,v_1)$ is an edge implies $\theta_1+\theta_4>0$ and the fact that $(v_1,v_2)$ is an edge implies $\theta_1+\theta_5>0.$ Similarly $v_1=(\theta_5,\{\theta_1,\theta_2\},\{\theta_2,\theta_3\}) $ implies $\theta_1+\theta_4>0$ and $\theta_1+\theta_5>0$ as well. Finally if $v_{-1}=(\theta_3,\{\theta_4,\theta_5\},\{\theta_1,\theta_2\}) $ the fact that $(v_{-1},v_0)$ is an edge implies $\theta_1+\theta_2>0$ and $\theta_1+\theta_3>0.$ Similarly if $v_{-1}=(\theta_2,\{\theta_4,\theta_5\},\{\theta_1,\theta_3\}) $  implies $\theta_1+\theta_2>0$ and $\theta_1+\theta_3>0$ and the lemma is proved. $\square$

\vspace{4mm}\noindent\textbf{Proof of Theorem 5.1.}  Suppose that $v_0$ is in a cycle  of $G_{\theta}^*.$ Therefore the Lemma 5.2 is applicable to $v_0.$ Recall now that in the statement of the theorem  we have assumed $\theta_2\leq \theta_3$ and $\theta_4\leq \theta_5.$ If $v_1=(\theta_4,\{\theta_1,\theta_5\},\{\theta_2,\theta_3\}) $ , since $(v_2,v_3)$ is an edge we get the new inequality $\theta_2+\theta_4>0.$ Similarly if $v_1=(\theta_5,\{\theta_1,\theta_2\},\{\theta_2,\theta_3\}) $ implies $\theta_2+\theta_5>0.$ Since $\theta_4\leq \theta_5$ the inequality $\theta_2+\theta_5>0$ holds for both possible values of $v_1.$

At this point observe that the only inequality to prove now is $\theta_3+\theta_4>0$, which imply the  other one $\theta_3+\theta_5>0.$
We discuss again the two possible values of $v_1$ for applying Lemma 5.2 to the sequence $ v_{-1}\rightarrow v_{0}\rightarrow v_{1}\rightarrow v_{2}\rightarrow v_3.$ If $v_1=(\theta_4,\{\theta_1,\theta_5\},\{\theta_2,\theta_3\}) $ we have seen that  $\theta_2+\theta_4>0$ which imply $\theta_3+\theta_4>0.$ 
Therefore the result is proved in this case. 

Now we assume  $v_1=(\theta_5,\{\theta_1,\theta_2\},\{\theta_2,\theta_3\}) $ and we discuss according to the two possible values of $v_2.$ If  $v_2=(\theta_2,\{\theta_3,\theta_5\},\{\theta_1,\theta_4\}) $ we apply Lemma 5.2 to the sequence $  v_{0}\rightarrow v_{1}\rightarrow v_{2}\rightarrow v_3\rightarrow v_4.$ Among the six inequalities we find 
$\theta_2+\theta_4>0$ which implies $\theta_3+\theta_4>0.$ The last case is $v_2=(\theta_3,\{\theta_2,\theta_5\},\{\theta_1,\theta_4\}) $: among the inequalities given by Lemma 5.2 we find the desired one $\theta_3+\theta_4>0.$ Finally from (\ref{change}) $BH_v$ always exists if $v$ is in a cycle of $G^*_{\theta}$ since in this case all the traditional parameters $a,b,p,q,r$ are positive.

One can observe  from the preceeding study that if $\theta_2+\theta_4\leq 0$ then necessarily $v_1=(\theta_5,\{\theta_1,\theta_2\},\{\theta_2,\theta_3\})$ and $v_2=(\theta_3,\{\theta_2,\theta_5\},\{\theta_1,\theta_4\}) .$ Iterating this remark we see that the set $(v_0,v_1,v_2,v_3,v_4,v_5)$ described in (\ref{VCC}) is a cycle.

We now discuss the converse. Suppose that $\theta_i+\theta_j>0$ for all $1\leq i<j\leq 5.$ In this case $G=G^*_{\theta}$ and since any vertex of $G$ is the initial vertex and the end vertex of some arrows, any vertex of $G$ belongs to a cycle.  Now suppose that $\theta_i+\theta_j>0$ for all $1\leq i<j\leq 5$ except for $(i,j)=(2,4).$ We have described above the corresponding cycle. For seeing that $\min \{\theta_1,\theta_3,\theta_5\}>\max\{\theta_2,\theta_4\}$
observe that we have seen that in an element of the cycle $\theta_i+\theta_j\leq 0$ can happen only for one pair $(i,j).$ For instance for $v_0$  we had $(i,j)=(2,4)$ and therefore $\theta_2<\theta _3$ and $\theta_4<\theta _5$ must be strict inequalities. Now  making the same remark for the other $(v_1,v_2,v_3,v_4,v_5)$ of (\ref{VCC}) shows $\min \{\theta_1,\theta_3,\theta_5\}>\max\{\theta_2,\theta_4\}.$ 
 $\square$

\vspace{4mm}\noindent \textsc{Cycles when $\theta_2+\theta_4\leq 0.$}
Assume as in Theorem 5.1 that $\theta_2\leq \theta_3$ and $\theta_4\leq \theta_5$ and that $\theta_2+\theta_4\leq 0.$ In this case $v_0$ as in Theorem 5.1 can belong to a cycle of size 1,2,3 or 6. Let us indicate here without proof the necessary and sufficient conditions for this. This is done by analyzing (\ref{VCC}).
\begin{itemize}
\item $v_0$ is in a cycle of size 1 if and only if  $v_0=x|xy|xy$ with $x+y>0$ and $y\leq 0.$
\item $v_0$ is in a cycle of size 2 if and only if  $v_0=x|xy|xz$ with $x+y, x+z>0,$  $y+z\leq 0$ and $y\neq z.$
\item $v_0$ is in a cycle of size 3 if and only if  $v_0=y|xz|xu$ with $x+y, x+z, x+u>0,$  $x\leq 0$ and $y,z,u$ are not all equal. 
\item $v_0$ is in a cycle of size 6 if and only if  $v_0=x|yz|uv$ with   $y+u\leq 0$, all the other sums of pairs are positive   $y\neq u$ and $x,z,v$ are not all equal. 
\end{itemize}

\vspace{4mm}\noindent \textsc{Example.} Let us use the identity $_2F_1(p,-p;1/2;\sin^2 \theta)=\cos 2p \theta.$ It provides   a hypergeometric function which is certainly positive on $(0,1)$ if $0<p\leq 1/2.$ Consider the distribution $BH(2\alpha,b,p,-p,1/2)$ for $\alpha, b>0.$ From (\ref{ITT}) it is equal to $BH_{\theta}$ with 
$$\theta_1=\alpha+b,\ \theta_2=-\alpha,\ \theta_3=-\alpha+\frac{1}{2},\ \theta_4=\alpha-p,\ \theta_5=\alpha+p$$
Note that $\theta_2<\theta_3$, $\theta_4<\theta_5$ and $\theta_2+\theta_4=-p<0.$
It is easy to detect with the help of Theorem 5.1 that that $v_{\theta}$ belongs to a cycle of $G^*_{\theta}$ if and only if $1/4> \alpha$
and $b+2\alpha> p,$ with $0<p<1/2.$ Since 
$\theta_2+\theta_4<0$ this cycle is described by (\ref{VCC}).

\vspace{4mm}\noindent \textsc{A cycle of order 6 changed in a cycle of order 5.} The above analysis includes $\theta_2+\theta_4=0.$ In this case $BH_{v_j}$ for $j=0,\ldots,5$ are beta distributions or quasi beta distributions. Since in this case the same distribution has an
infinite number of $\theta$ parameters, it makes sense
to ask if these different representations could appear within a
cycle. This happens only in the following case.  Taking $v_0=(x,y, x+2y,-y,x)$ with $x>y>0$ we get $v_1=(x,-y,x,y,x+2y),$ $v_2=(x+2y,y,x,-y,x),$ $ v_3=(x,-y,x+2y,y,x),$ $v_4=(x,y,x,-y,x+2y),$   $v_5=(x+2y,-y,x,y,x)$  and 
$$BH_{v_0}=\beta_{x-y,x+y},\ BH_{v_1}=\beta_{x+3y,x-y},\ BH_{v_2}=\beta_{x-y,x+3y},$$$$\ BH_{v_3}=\beta_{x+y,x-y},\ BH_{v_4}=\beta_{x+y,x+y},\ BH_{v_5}=\beta_{x+y,x+y}.$$ We observe that $BH_{v_4}= BH_{v_5}.$ All the other $BH_{v_j}$ are different. Therefore, starting with a cycle of order 6 we can design a beta distributed random continued fraction of the type of Theorem 4.3 with period $k=5.$ One can prove that this phenomenon appears only for the above choice of parameters. 

\section{Aknowledgments} We thank   Sapienza Universit\`a di Roma   and the Universit\'e Paul Sabatier in Toulouse for their generous support during the preparation of this paper. We are grateful to  Giovanni Sebastiani of IAC-CNR in Rome who has done the computer calculations for the cycles of the graph 1+1+1+1+1. We are also grateful to  Richard Askey who pointed us to  the work of Felix Klein and to  Jean-Fran\c{c}ois Chamayou who introduced us to the Thomae formula and helped us read the paper by Klein.

\section{References}\vspace{4mm}\noindent\textsc{Abramovitz, M. and Stegun, I. A.} (1965) \textit{Handbook of Mathematical Functions}, Dover Publications, New York.

 \vspace{4mm}\noindent \textsc{Andrews, G. E.,  Askey, R. and Roy, R.} (1999) \textit{Special Functions}, Encyclopedia Math. Appl., vol. 71, Cambridge Univ. Press, Cambridge. 

\vspace{4mm}\noindent\textsc{Asci, C., Letac, G. and Piccioni, M.} (2008)
'Beta-hypergeometric distributions and random continued fractions', 
\textit{Statist. Probab. Lett.} \textbf{78}, 1711-1721.

 \vspace{4mm}\noindent\textsc{Bailey, W. N.} (1935) \textit{Generalized Hypergeometric Series}, Cambridge Mathematical Tract no. 32, Cambridge. 

\vspace{4mm}\noindent \textsc{Beyer, W.A., Louck, J. D. and Stein, P. R.} (1987) 'Group
theoretical basis of some identities for the generalized hypergeometric
series', \textit{J. Math. Phys. }\textbf{28,} 497-508.

 \vspace{4mm}\noindent \textsc{Chamayou, J.-F. and Letac, G.} (1991) 'Explicit stationary
distributions for composition of random functions and products of random
matrices', \textit{J. Theor. Probab.} \textbf{4, }3-36\textbf{.}

\vspace{4mm}\noindent \textsc{Chrystal, G.} (1964) \textit{Textbook of Algebra,} seventh edition, vol. 2 Chelsea, New York. 

 \vspace{4mm}\noindent \textsc{Dyson, F. J.} (1953) 'The dynamics of a disordered linear chain', \textit{ Phys. Rev.}  \textbf{92}, 1331-1338.

\vspace{4mm}\noindent \textsc{Hurwitz, A.} (1891) 'Ueber die Nullstellen der hypergeometrishen Reihe', \textit{Math. Ann.} \textbf{38}, 452-458.

\vspace{4mm}\noindent \textsc{Klein, F.} (1890) 'Ueber die Nullstellen der hypergeometrishen Reihe', \textit{Math. Ann.} \textbf{37}, 573-590.

\vspace{4mm}\noindent \textsc{Letac, G.} (1986)
 'A contraction principle for certain Markov chains and its applications',
 {\it Contemp. Math.} {\bf 50}, 263-273.

 \vspace{4mm}\noindent \textsc{Letac, G. and  Seshadri, V.} (1983)
 'A characterization of the generalized inverse Gaussian distribution by continued fractions',
 {\it Z. Wahrsheinlichkeitstheorie und Verv. Geb. } {\bf
62}, 485-489.

 \vspace{4mm}\noindent \textsc{Lloyd, P.} (1969) 'Exactly solvable model in electronic states in a three dimensional disordered Hamiltonian: non-existence of localized states', \textit{J. Phys. C. Solid State Phys.} \textbf{2} 1717-1725.

 \vspace{4mm}\noindent \textsc{Maier, R. S.} (2005) 'A generalization of Euler's
transformation', \textit{Trans. Amer. Math. Soc}. \textbf{358}, 39-57.

\vspace{4mm}\noindent \textsc{Marklov, J.,  Tourigny, Y. and  Wolowski, L.} (2008)) 'Explicit invariant measures for products of random matrices', \textit{Trans. Amer. Math Soc.} \textbf{360} (2008) 3391-3427.

\vspace{4mm}\noindent \textsc{Propp, J. and Wilson, D.} (1996) 'Exact Sampling with Coupled Markov Chains and Applications to Statistical Mechanics,' \textit{Random Structures Algorithms } \textbf{9} 223-252.

 \vspace{4mm}\noindent \textsc{Rainville, E. D.} (1960) \textit{Special Functions}, Macmillan, New
York.

\vspace{4mm}\noindent \textsc{Thomae, J.} (1879) 'Ueber die Funktionen, welche durch
Reichen von der Form dargestellst verden $1+\frac{p}{1}\frac{p^{\prime }}{%
q^{\prime }}\frac{p^{\prime \prime }}{q^{\prime \prime }}+\frac{p}{1}\frac{%
p+1}{2}\frac{p^{\prime }}{q^{\prime }}\frac{p^{\prime }+1}{q^{\prime }+1}%
\frac{p^{\prime \prime }}{q^{\prime \prime }}\frac{p^{\prime \prime }+1}{%
q^{\prime \prime }+1}+\cdots ,$' \textit{J. Reine Angew. Math.} \textbf{87} ,
26-73.

\vspace{4mm}\noindent \textsc{Titchmarsh, E. C.} (1939) \textit{Theory of Functions}, second edition. Oxford University Press, Amen House, London.

\vspace{4mm}\noindent \textsc{Van Vleck, E. B.} (1902) 'A determination of the number of real and imaginary roots of the hypergeometric series', \textit{Trans. Amer. Math. Soc. } \textbf{3}, 110-131.

\end{document}